\input amstex
\documentstyle{amsppt}
\magnification=\magstep1 \NoRunningHeads
\topmatter
\title
Actions of finite  rank:
weak rational ergodicity and partial rigidity
\endtitle
\author
Alexandre I. Danilenko
\endauthor

\address
 Institute for Low Temperature Physics
\& Engineering of National Academy of Sciences of Ukraine, 47 Lenin Ave.,
 Kharkov, 61164, UKRAINE
\endaddress
\email alexandre.danilenko\@gmail.com
\endemail

\abstract
A simple proof of the fact that each rank-one infinite measure preserving (i.m.p.) transformation is subsequence weakly rationally ergodic is found.
Some classes of funny rank-one i.m.p. actions of Abelian groups are shown to be subsequence weakly rationally ergodic.
A constructive definition of finite funny rank for actions of arbitrary infinite countable groups is given.
It is shown that the ergodic i.m.p. transformations of balanced finite funny rank are subsequence weakly rationally ergodic.
It is shown that the ergodic probability preserving transformations of exact finite rank, the ergodic Bratteli-Vershik maps corresponding to the ``consequtively ordered'' Bratteli diagrams of finite rank, some their generalizations and the ergodic IETs are partially rigid.
\endabstract

 \loadbold
\endtopmatter

\document

\head 0. Introduction
\endhead

In \cite{Aa1},  Aaronson (motivated by estimating the asymptotic growth of  some special ergodic averages in infinite measure spaces) introduced  concepts of {\it rational} and  {\it weak rational ergodicity} for infinite measure preserving (i.m.p.) transformations.
More recently, in \cite{Aa2}, he introduced a related concept of   rational weak mixing and considered more general {\it subsequence} versions of these concepts.
All the examples of systems possessing these properties given in \cite{Aa1} and \cite{Aa2} 
have positive Krengel   entropy and  countable Lebesgue spectrum.
Other kind of examples---rank-one i.m.p. transformations---related to  the aforementioned  concepts were considered in
 subsequent papers \cite{Aa3},  \cite{Dai--Si}  and \cite{Bo--Wa}.
We recall that the rank-one transformations have zero Krengel entropy and simple spectrum  (see \cite{DaSi} and references therein).

The main result of \cite{Bo--Wa} is that each rank-one i.m.p. transformation is subsequence weakly rationally ergodic.
In the present paper we give a short proof of this fact.
Moreover, we generalize  the concept of subsequence weak rational ergodicity to i.m.p. actions of arbitrary countable discrete groups and construct some families  of weakly rationally ergodic  (along a F{\o}lner sequence) funny rank-one i.m.p. actions for  arbitrary countable infinite Abelian groups.
We obtain as an immediate corollary that all these actions  are  non-squashable, i.e. every nonsingular transformation commuting with this action preserves the measure.
This fact  (i.e. the non-squashability)  was  first proved in \cite{Aa1} for  the weakly rationally ergodic i.m.p. $\Bbb Z$-actions  but  that proof is valid in the general case of weakly rationally ergodic  (along a F{\o}lner sequence) i.m.p. group actions.

One of the main result of this paper is that the ergodic i.m.p. transformations of balanced finite rank are subsequence weakly rationally ergodic.
We say that an i.m.p.  transformation  is of {\it balanced} finite rank if it is of finite rank and the bases of the Rokhlin towers  on the $n$-th step of the cutting-and-stacking inductive construction have  asymptotically comparable measures as $n\to\infty$.  
We recall that the transformations of finite rank have zero Krengel entropy and finite spectral multiplicities (see \cite{DaSi}).

We give a constructive definition of  finite funny rank actions for arbitrary  countable infinite groups by developing the $(C,F)$-construction  defined originally in \cite{dJ} and \cite{Da1} to produce funny rank-one actions.
We hope that this generalized $(C,F)$-construction will find other applications in ergodic theory, especially in the spectral theory, the theory of joinings of dynamical systems, the theory of i.m.p. and nonsingular systems, etc.  (see the survey \cite{Da3} for various application of the $(C,F)$-techniques in the rank-one case).

Rosenthal  showed in an unpublished paper  \cite{Ro} that the ergodic (finite measure preserving) transformation of exact finite rank is not mixing. 
Recently this fact was reproved under some restriction in \cite{Be--So}.
The {\it exactness} means that the transformation is constructed via the cutting-and-stacking procedure without adding
spacers and the corresponding Rokhlin towers do not asymptotically vanish.
Since the exact rank-one transformations have pure point spectrum,
 the transformations of exact rank greater  than one can be thought of as a ``higher rank'' analogues of systems with pure point spectrum. 
We refine Rosenthal's result by showing that the ergodic  transformations of exact finite rank are partially rigid\footnote{When this paper had been already submitted, V.~Ryzhikov informed the author about his earlier work \cite{Ry}. Though the partial rigidity of the finite rank transformations was not explicitly asserted there, it was actually proved there.}.
We also extend this assertion to the transformations of quasi-exact finite rank, which means that spacers in the underlying cutting-and-stacking construction are possible but with some uniform (over the indicative steps) bound on their number. 
It was proved in a recent paper \cite{Be--So} that some Vershik transformations associated with the so-called {\it consecutively ordered} (see \cite{Du} for the definition) Bratteli diagrams of finite rank  are non-mixing.
We show that these transformations (and some generalizations of them) are indeed partially rigid.
It is also shown how to deduce from  Katok's proof \cite{Ka} of non-mixing for the ergodic interval exchange transformations (IETs) that they are partially rigid\footnote{The fact that the ergodic IETs are partially rigid was also established (implicitly) in \cite{Ry}.}.

The paper is organized as follows.
Section~1 is devoted to the $(C,F)$-construction of funny rank-one actions for discrete countable groups.
The construction appeared first in \cite{dJ} and \cite{Da1} in slightly different versions (see also  \cite{Da3}). 
We compare them in the present paper and introduce a third version which is formally more general (in fact, the most possible general in view of  Proposition~1.4) than the two ones.
We show however in Theorem~1.8 that the class of measurable $(C,F)$-actions in the sense of the third definition (which is exactly the class of all  funny rank-one actions with an invariant $\sigma$-finite measure by Theorem~1.6) is the same as the class of measurable $(C,F)$-actions in the sense of the definition from \cite{Da1} if the actions are considered up to modification on null subsets.
Moreover, in the {\it finite} measure preserving case the three versions of the $(C,F)$-constructions define the very same class of actions  and the acting groups are necessarily amenable.
Since every $(C,F)$-action $T$ in the sense of the definition from \cite{Da1} is {\it strictly ergodic}, i.e.  $T$ is a topological action on a locally compact second countable space, $T$ is minimal and  $T$ admits a unique  up to scaling invariant $\sigma$-finite Radon measure, we obtain as a byproduct strictly ergodic models for arbitrary funny rank-one  $\sigma$-finite measure preserving actions of arbitrary discrete countable groups.
We recall that strictly ergodic models for the arbitrary (not only rank-one) ergodic finite measure preserving actions of Abelian groups were constructed in 
\cite{We1} (see also a discussion there for earlier results) and  strictly ergodic models for the ergodic i.m.p.  $\Bbb Z$-actions were constructed in \cite{Yu2}. 
We also mention another application of the $(C,F)$-construction.
By \cite{Zi}, if a discrete countable group $G$ admits a free ergodic probability preserving action whose orbit equivalence relation  is hyperfinite then $G$ is amenable.
However,  each non-amenable group has  i.m.p.  free actions with hyperfinite orbit equivalence relations (see e.g., \cite{BeGo}).
 The $(C,F)$-construction provides a simple  way to obtain such actions  possessing    additional properties such as a strict ergodicity (in locally compact spaces), funny rank one, etc.

In Section~2, for an arbitrary discrete countable amenable group $G$, we introduce a concept of weak rational ergodicity along a F{\o}lner sequence in $G$.
In the case $G=\Bbb Z$ and the F{\o}lner sequence consists of  intervals  with $0$ as the left endpoint, our concept coincides with Aaronson's sequence  weak rational ergodicity \cite{Aa2}. 
Using the language of the $(C,F)$-construction we give a short proof for the main result of the first version of  \cite{Bo--Wa} that each rank-one i.m.p. transformation is subsequence weakly rationally ergodic (Theorem~2.4)\footnote{Being informed about our proof of Theorem~2.4, the authors of \cite{Bo-Wa} replaced their main result with a stronger one in the final version of \cite{Bo-Wa}. The two versions of  \cite{Bo-Wa} can be found in ArXiv.} .
Then we extend this result to some classes of funny rank-one i.m.p. actions of arbitrary  countable Abelian groups (Corollary~2.7 and Theorem~2.8).
The question whether every funny rank-one i.m.p. transformation (or an Abelian group action) is
subsequence weakly rationally ergodic remains open.

In Section 3 we consider  $\sigma$-finite measure preserving group actions of finite funny rank 
 (see \cite{Fe2} for the definition of finite funny rank in the case when the acting group is $\Bbb Z$).
 We pay special attention to ergodic $\Bbb Z$-actions of finite  rank.
 It is shown that each ergodic $\sigma$-finite measure preserving transformation of finite rank is  built over a finite measure preserving  transformation of exact finite rank and under a piecewise constant integer valued function (see Corollary~3.5).

In Section~4 we generalize the $(C,F)$-construction in such a way that it yields  actions of finite funny rank and that every $\sigma$-finite measure preserving action of  finite funny rank is isomorphic to a $(C,F)$-action  in the generalized sense (see Theorems~4.8 and 4.9). 
Thus the $(C,F)$-construction can be considered as a constructive definition for the actions of finite funny rank.
As in the rank-one case, we obtain strictly ergodic models for the ergodic $\sigma$-finite measure preserving  $G$-actions of finite funny rank (Theorem~4.10). 
In the case where $G=\Bbb Z$ and the actions are of finite rank, we show how to associate an ordered Bratteli diagram to  a $(C,F)$-data in such a way that the corresponding $(C,F)$-action of $G$ is the Bratteli-Vershik map associated with the diagram (Remark~4.11).
Thus the $(C,F)$-construction can be viewed as a generalization of Bratteli-Vershik construction from $\Bbb Z$-actions to actions of arbitrary discrete countable groups.\footnote{In this paper we consider  only $(C,F)$-actions of finite funny rank. However the infinite funny rank $(C,F)$-actions can be defined in a similar way. Remark~4.11 will hold for them as well. This will be done elsewhere.}

In Section~5  we show that  the ergodic i.m.p. transformations of  balanced finite rank are weakly rationally ergodic (Theorem~5.5).

In Section~6 we consider  finite measure preserving transformations of finite rank.
We show that each ergodic transformation of exact finite rank is partially rigid (Theorem~6.1).
Then in Theorem~6.4 we extend it to the transformations of quasi-exact finite rank. 
In a similar way we  strengthen another result from \cite{Be--So} on non-mixing for another class of transformations of finite rank.
We prove in Theorem~6.6 that each ergodic transformations of finite rank  with  consecutive ordering of towers that satisfies  a ``non-degeneracy'' condition is partially rigid.
It is also shown that the ergodic IETs are partially rigid (Proposition~6.8).

The final Section~7 is a list of open problems related to the subject of this paper.

{\sl Acknowledgements.} The author thanks M.~Lema{\'n}chyk  and K.~Fr{\c a}czek for useful discussions.


\head 1.   $(C,F)$-construction of funny rank-one  actions 
\endhead

Let $T=(T_g)_{g\in G}$ be a measure preserving  action of  a countable infinite discrete group $G$ on a standard $\sigma$-finite measure space $(X,\goth B,\mu)$.
The following definition was given by J.-P.~Thouvenot in the case $G=\Bbb Z$ (see also \cite{Fe1} and \cite{Fe2}).

\definition{Definition 1.1}
If there exist a sequence  $(B_n)_{n\ge 0}$ of subsets of finite measure  in $X$ and a sequence  $(F_n)_{n\ge 0}$ of finite subsets in $G$ such that  
\roster
\item"(i)" 
for each $n\ge 0$, the subsets $T_gB_n$, $g\in F_n$, are pairwise disjoint and 
\item"(ii)" 
for each subset $B\in\goth B$ with $\mu(B)<\infty$,
$$
\lim_{n\to\infty}\inf_{F\subset F_n}\mu\bigg(B\triangle\bigsqcup_{g\in F}T_gB_n\bigg)=0
$$
\endroster
then $T$ is called  an action {\it of funny rank one}.
\enddefinition

If $G=\Bbb Z$ and every $F_n$ is an interval $\{0,1,\dots, \#F_n-1\}$ in $\Bbb Z$ then $T$ is called an action {\it of rank one}.
For the constructive definition of rank-one transformations (i.e. $\Bbb Z$-actions) using the cutting-and-stacking inductive process we refer to \cite{Fe2}.
A constructive definition of actions of funny rank-one  was given in \cite{dJ} and \cite{Da1} (see also  \cite{Da2}).
We now recall it.

Let $(F_n)_{n\ge 0}$  and $(C_n)_{n\ge 1}$ be two sequences of finite subsets in $G$ such that for each $n>0$,
\roster
\item"(I)"
$F_0=\{1\}$,  $\# C_n>1$,
\item"(II)"
$F_nC_{n+1}\subset F_{n+1}$,
\item"(III)"
$F_nc\cap F_nc'=\emptyset$ if $c,c'\in C_{n+1}$ and $c\ne c'$.
\endroster
We let $X_n:=F_n\times C_{n+1}\times C_{n+2}\times\cdots$ and endow this set with the infinite product topology.
Then $X_n$ is a compact Cantor (i.e. totally disconnected  perfect metric) space.
The  mapping
$$
X_n\ni (f_n,c_{n+1},c_{n+2},\dots)\mapsto(f_nc_{n+1}, c_{n+2},\dots)\in X_{n+1}
$$
is a topological embedding
of $X_n$ into $X_{n+1}$.
Therefore an inductive limit $X$ of the sequence $(X_n)_{n\ge 0}$ furnished with these embeddings  is a well defined locally compact Cantor  space.
We call it the {\it $(C,F)$-space associated with the sequence} $(C_n,F_{n-1})_{n\ge 1}$.
It is easy to see that the $(C,F)$-space is compact if and only if there is $N>0$ with $F_{n+1}=F_nC_{n+1}$ for all $n>N$.
Given a subset $A\subset F_n$, we let 
$$
[A]_n:=\{x=(f_n,c_{n+1},\dots)\in X_n\mid f_n\in A\}
$$ 
and call this set an  {\it $n$-cylinder} in $X$.
It is open and compact in $X$.
The collection of all cylinders coincides with the family of  all compact open subset  in $X$.
It is easy to see that
$$
\align
[A]_n\cap[B]_n &=[A\cap B]_n, \quad
[A]_n\cup[B]_n =[A\cup B]_n\quad \text{and}\\
 [A]_n &=[AC_{n+1}]_{n+1}
 \endalign
 $$
  for all $A,B\subset F_n$ and  $n\ge 0$.
For brevity, we will write $[f]_n$ for $[\{f\}]_n$, $f\in F_n$.

Let $\Cal R$ denote the {\it tail equivalence relation} on $X$.
This means that the restriction of $\Cal R$ to $X_n$ is the tail equivalence relation on $X_n$ for each $n\ge 0$.
We note that $\Cal R$ is {\it minimal}, i.e. the $\Cal R$-class of every point  is dense in $X$.
There exists a unique $\sigma$-finite $\Cal R$-invariant Borel measure $\mu$ on $X$ such that
$\mu(X_0)=1$.
It is a Radon measure, i.e. it is finite on every compact subset.
Moreover,  $\mu$ is strictly positive on every non-empty open subset.
We note that the {\it $\Cal R$-invariance}\footnote{$\mu$ is called $\Cal R$-invariant if $\mu$ is invariant under each Borel transformation whose graph is contained in $\Cal R$.} of $\mu$ is equivalent to the following property:
$$
\mu([f]_n)=\mu([f']_n)\quad\text{ for all }f,f'\in F_n, n\ge 0.
$$
It is easy to see that 
$$
\mu([A]_n)=\frac{\# A}{\# C_1\cdots\# C_n} \quad\text{for each subset }A\subset F_n, n>0.
$$
We call $\mu$ the {\it $(C,F)$-measure associated with $(C_n, F_{n-1})_{n\ge 1}$}.
It is finite if and only if\footnote{In view of (I)--(III), the sequence $(\frac{\# F_n}{\# C_1\cdots\#C_n})_{n=1}^\infty$ is non-decreasing and bounded by $1$ from below.} 
$$
\lim_{n\to\infty}\frac{\# F_n}{\# C_1\cdots\#C_n}<\infty.
\tag1-1
$$
It is easy to see that  $\mu$  on is {\it $\Cal R$-ergodic}, i.e. each Borel $\Cal R$-saturated subset of $X$ is either $\mu$-null or $\mu$-conull.
We now define an action of $G$ on $X$ (or, more rigorously,  on a subset of $X$).
Given $g\in G$, let
$$
X_n^g:=\{(f_n,c_{n+1},c_{n+2}\dots)\in X_n\mid gf_n\in F_n\}.
$$
Then $X_n^g$ is a compact open subset of $X_n$ and $X_n^g\subset X_{n+1}^g$.
Hence the union $X^g:=\bigcup_{n\ge 0}X_n^g$ is  a well defined  open subset of $X$.
Let $X^G:=\bigcap_{g\in G}X^g$.
Then $X^G$ is a $G_\delta$-subset of $X$.
Hence $X^G$ is Polish in the induced topology.
Given $x\in X^G$ and $g\in G$, there is $n>0$ such that $x=(f_n,c_{n+1},\dots)\in X_n$ and $g f_n\in F_n$.
We now let 
$$
T_gx:=(gf_n, c_{n+1},\dots)\in X_n\subset X.
$$
It is standard to verify that 
\roster
\item"(i)"$T_gx\in X^G$,
\item"(ii)" the map $T_g:X^G\ni x\mapsto T_gx\in X^G$ is a homeomorphism of $X^G$ and 
\item"(iii)"$T_gT_{g'}=T_{gg'}$ for all $g,g'\in G$.
\endroster
Thus $T:=(T_g)_{g\in G}$ is a continuous action of $G$ on $X^G$.
We call it {\it the $(C,F)$-action of $G$ associated with the sequence $(C_n,F_{n-1})_{n\ge 0}$}.
It is free.
It is obvious that $X^G$ is $\Cal R$-invariant and the $T$-orbit equivalence relation is the restriction of $\Cal R$ to $X^G$.
It follows that  $T$ preserves  $\mu$.

Given sequences $(F_n)_{n\ge 0}$ and $(C_n)_{n\ge 1}$ satisfying (I)--(III) and a sequence $(z_n)_{n\ge 1}$ of elements of $G$, we let
$C_n':=z_n^{-1}C_nz_{n+1}^{-1}$ and $F_{n-1}':=F_{n-1}z_{n}$ for each $n\ge 1$.
Then the sequences $(F_n')_{n\ge 0}$ and $(C_n')_{n\ge 1}$ satisfy (I)--(III).
Denote by $X'$, $\Cal R'$ and $T'$ the associated $(C,F)$-space,  $(C,F)$-equivalence relation and   $(C,F)$-action respectively.
Then there is a canonical homeomorphism $\phi:X\to X'$ that intertwines $\Cal R$ with $\Cal R'$ and $T$ with $T'$.
It is given by
$$
\phi(f_n,c_{n+1},c_{n+2},\dots)=(f_nz_{n+1}, z_{n+1}^{-1}c_{n+1}z_{n+2}, z_{n+1}^{-1}c_{n+1}z_{n+2},\dots )\in X_n'\subset X'
$$
whenever $(f_n,c_{n+1},c_{n+2},\dots)\in X_n\subset X $, for each $n\ge 0$.
Choosing $(z_n)_{n=1}^\infty$ in an appropriate way we may always assume without loss of generality\footnote{This means that we can modify the $(C,F)$-sequences in such a way that   the modified associated $(C,F)$-action (and the $(C,F)$-equivalence relation) is isomorphic to the original one.} that the following condition
\roster
\item"(IV)" $1\in\bigcap_{n\ge 0}F_n\cap\bigcap_{n\ge 1}C_n$
\endroster
is always satisfied in addition to (I)--(III).

 Since $X^G$  is $\Cal R$-saturated  and  $\mu$ is $\Cal R$-ergodic, we have
either $\mu(X^G)=0$ or  $\mu(X\setminus X^G)=0$.
Each of the two cases is possible to occur.
Moreover,   $X^G$ can be empty at all.
 We now discuss conditions under which   $X^G$ is $\mu$-conull  or even $X^G=X$.  

\proclaim{Proposition 1.2}
 $X^G=X$ if and only if for each $g\in G$ and $n>0$, there is $m>n$ such that
$$
gF_nC_{n+1}C_{n+2}\cdots C_{m}\subset F_m.\tag1-2
$$
\endproclaim
\demo{Proof}
The ``if'' part is trivial.
We now  prove the ``only if".
Fix $n\ge 0$ and $g\in G$.
The map 
$$
\phi_m:X_n\ni x=(f_n,c_{n+1},\dots)\mapsto gf_nc_{n+1}\cdots c_m\in G
$$
is continuous, $m>n$.
Since $X^g=X$, we obtain that $X_n=\bigcup_{m> n}\phi_m^{-1}(F_m)$.
Since $X_n$ is compact and $\phi_{n+1}^{-1}(F_{n+1})\subset \phi_{n+2}^{-1}(F_{n+2})\subset\cdots$, 
it follows that $X_n=\phi_{m}^{-1}(F_{m})$ for some $m>n$.
The inclusion \thetag{1-2} follows.
\qed
\enddemo

Thus, in this case the $(C,F)$-action is defined on the entire (locally compact) space $X$.

\remark{Remark \rom{1.3}}
We note that if $X$ is not compact then  $T$ extends to the one-point compactification $X^*=X\sqcup\{\infty\}$  of $X$ by setting $T_g\infty=\infty$ for all $g\in G$.
We obtain a continuous action of $G$ on $X^*$.
This action is {\it almost minimal}, i.e. there is one fixed point and the orbit of any other point is dense.
This concept was introduced in \cite{Da2} in the case $G=\Bbb Z$. 
For the (topological) orbit classification of the almost minimal $\Bbb Z$-systems see \cite{Da2} and \cite{Ma}.
Some natural examples of such systems (subshifts arising from non-primitive substitutions) are given in \cite{Yu1}.
\endremark

\proclaim{Proposition 1.4}
 $\mu(X\setminus X^G)=0$
 if and only if for each $g\in G$ and $n>0$,
$$
\lim_{m\to\infty}\frac{\#((gF_nC_{n+1}C_{n+2}\cdots C_{m})\cap F_m)}{\# F_n\#C_{n+1}\cdots\# C_m}=1.
\tag1-3
$$
If $\mu(X)<\infty$ then $\mu(X\setminus X^G)=0$ if and only if 
$(F_n)_{n\ge 0}$ is a F{\o}lner sequence in $G$ and hence $G$ is amenable\footnote{Another way to see that $G$ is amenable is to apply a theorem by R. Zimmer from \cite{Zi}: if  $G$ has a free probability preserving ergodic action $T$ and the $T$-orbit equivalence relation is hyperfinite then $G$ is amenable. Of course, the tail equivalence relation is hyperfinite, i.e. it is the union of an increasing sequence of equivalence relations with finite equivalence classes.}.
\endproclaim

\demo{Proof}
Since $\mu(X\setminus X^G)=0$ if and only if $\frac{\mu(X_n\cap X_m^g)}{\mu(X_n)}\to 1$ as $m\to\infty$ for each $g\in G$ and $n>0$, it suffices to note that 
$$
\frac{\mu(X_n\cap X_m^g)}{\mu(X_n)}=\frac{\#((gF_nC_{n+1}\cdots C_{m})\cap F_m)}{\# F_n\#C_{n+1}\cdots\# C_m}.
$$
 In the case where $\mu$ is finite we have $\mu(X\setminus X^G)=0$ if and only if
 $\frac{\mu(X_n^g)}{\mu(X_n)}\to 1$ as $n\to\infty$ for each $g\in G$.
 Since
 $$
 \frac{\mu(X_n^g)}{\mu(X_n)}=\frac{\#(gF_n\cap F_n)}{\# F_n},
 $$
 the second assertion of the lemma follows. \qed 
\enddemo

We note that the condition $\mu(X)<\infty$ from the second claim of Proposition~1.4 can not be omitted. Indeed, if  $G$ is amenable then it is not difficult to construct
 sequences $(C_n,F_{n-1})_{n\ge 1}$ such that $(F_n)_{n\ge 0}$ is F{\o}lner, (I)--(IV)  are satisfied  but~\thetag{1-3}  is not satisfied.
Hence $\mu(X^G)=0$ and
  $\mu(X)=\infty$.

From now on we assume that \thetag{1-3} holds.
Then $(X,\mu,T)$ is a $\sigma$-finite measure preserving dynamical system.
Of course, it is free, conservative and ergodic.
We claim that it is of funny rank one.
Indeed, the  sequences $([1]_n)_{n\ge 0}$ and $(F_n)_{n\ge 0}$ 
satisfy Definition~1.1.
We note that 
$$
T_g[f]_n=[gf]_n \quad\text{(up to $\mu$-null subset) whenever $f,gf\in F_n$}.\tag1-4
$$

We summarize  the aforementioned  results on  $(C,F)$-actions  in the following theorem.

\proclaim{Theorem 1.5}
Given a sequence $(C_n,F_{n-1})_{n\ge 1}$ satisfying  (I)--(IV), there is a locally compact Cantor space $X$ and a countable equivalence relation $\Cal R$ on $X$ such that
\roster
\item"(i)"
every $\Cal R$-class is dense in $X$,
\item"(ii)"
there is only one (up to scaling) $\Cal R$-invariant non-trivial $\sigma$-finite Radon measure $\mu$ on $X$,
\item"(iii)"
$\mu$ is finite if and only if \thetag{1-1} is satisfied,
\item"(iv)" 
there is a free topological  $G$-action $T$ on an $\Cal R$-invariant $G_\delta$-subset $X^G$ of $X$ such that the $T$-orbit equivalence relation is the restriction of $\Cal R$ to $X^G$,
\item"(v)"
$X^G=X$ if and only if \thetag{1-2} is satisfied,
\item"(vi)"
$\mu(X\setminus X^G)=0$ if and only if \thetag{1-3} is satisfied.
Moreover, if  \thetag{1-3} is not satisfied then $\mu(X^G)=0$.
\item"(vii)"
If $\mu(X)<\infty$ then \thetag{1-3} is equivalent to the fact that $(F_n)_{n\ge 0}$ is a F{\o}lner sequence in $G$.
\item"(viii)"
Under  \thetag{1-3}, the dynamical system $(X,\mu,T)$ is ergodic, conservative and of funny rank one.
\endroster
\endproclaim

The converse to Theorem~1.5(viii) also holds.

\proclaim{Theorem 1.6 ($(C,F)$-models for actions of funny rank one)} If $T$ is a $\sigma$-finite measure preserving $G$-action of funny rank one then $T$ is isomorphic to a $(C,F)$-action of $G$ equipped with the $(C,F)$-measure.
\endproclaim
\demo{Proof}
Let $(B_n)_{n\ge 0}$ and $(F_n)_{n\ge 0}$ be as in Definition~1.1.
Without loss  of generality we may assume that  $F_0=\{1\}$, $1\in\bigcap_{n\ge 0}F_n$ and for each $g\in F_n$, the subset  $T_gB_{n}$ is  a union of subsets  $T_sB_{n+1}$   for some family of $s\in F_{n+1}$.
We now define inductively two sequences $(C_n)_{n\ge 1}$ and $(\widetilde F_n)_{n\ge 0}$ of finite subsets in $G$ that satisfy (I)--(III) and 
for each $n$,  
 $$
 \{T_gB_n\mid g\in F_n\}=\{T_gB_n\mid g\in \widetilde F_n\}.\tag1-6
 $$
For that, we first set $\widetilde F_0:=F_0$.
Suppose now that $\widetilde F_n$ is defined and \thetag{1-6} holds for some $n$.
There is a finite subset $C_{n+1}\subset F_{n+1}$ such that $B_n=\bigsqcup_{c\in C_{n+1}}T_cB_{n+1}$.
Then, in view of~\thetag{1-6}, and the refining property of towers, for each $g\in \widetilde F_n$, there is a finite subset $I_g\subset F_{n+1}$ such that
$$
\bigsqcup_{c\in C_{n+1}}T_{gc}B_{n+1}=T_gB_n=\bigsqcup_{c\in I_g}T_{g}B_{n+1}.
$$
Of course, the subsets $I_g$, $g\in\widetilde F_n$, are pairwise disjoint.
We now set $\widetilde F_{n+1}:=\widetilde F_nC_n\sqcup (F_{n+1}\setminus\bigsqcup_{g\in\widetilde F_n} I_g)$.
Then, of course, \thetag{1-6} is satisfied for $n+1$.
Moreover, 
it is easy to verify that the sequence $(C_n, \widetilde F_{n-1})_{n>1}$ satisfies (I)--(III).
We denote the corresponding $(C,F)$-action of $G$ by $S$.
In view of Definition~1.1(ii),
the one-to-one correspondence
$$
\bigsqcup_{g\in F} T_gB_n \mapsto [F]_n,\quad F\subset \widetilde F_n,
$$
between 
\roster
\item"(a)" the subsets measurable with respect to some of the partition $\{T_gB_n\mid g\in F_n\}$, $n>0$, and
\item"
(b)" the compact open subsets of the underlying $(C,F)$-space
\endroster
generates a Borel isomorphism (mod 0) between the underlying measure spaces.
Now \thetag{1-4} yields that
this isomorphism  intertwines $T$ with $S$.
Therefore $S$ is defined almost everywhere and hence \thetag{1-3} holds by Proposition~1.4.
\qed
\enddemo

We obtain the following corollary from Theorem~1.6 and Proposition~1.4.

\proclaim{Corollary 1.7}
If $G$ has a finite measure preserving action of  funny rank one then $G$ is amenable.
\endproclaim

It may seem that the condition \thetag{1-3} is essentially more general   than \thetag{1-2}.
However we  show that the two conditions determine    the same class of $(C,F)$-actions (modulo measure theoretic isomorphism).
First of all we introduce a technique of passing to a $(C,F)$-subsequence.
Let $T$ be the  $(C,F)$-action associated with  a sequence $(C_n,F_{n-1})_{n>0}$.
Given an increasing  sequence $(k_n)_{n\ge0}$ of non-negative integers with $k_0=0$, we let
$\widetilde F_n:=F_{k_n}$ and $\widetilde C_{n}:=C_{k_{n-1}+1}C_{k_{n-1}+2}\cdots C_{k_{n}}$.
 Since
$( C_n, F_{n-1})_{n>1}$ satisfies (I)--(III) and~\thetag{1-3},  the sequence $(\widetilde C_n,\widetilde F_{n-1})_{n>1}$  also satisfies these conditions.
We call the latter sequence a {\it $(C,F)$-subsequence} of $( C_n, F_{n-1})_{n>1}$.
Denote by  $\widetilde T$ 
 the $(C,F)$-action associated with it.
  Then $\widetilde T$ is canonically isomorphic to $T$.
  Indeed, let $X$ and $\widetilde X$ denote the corresponding $(C,F)$-spaces.
  We recall that 
  $$
  X=\bigcup_{n\ge 0}X_{n}=\bigcup_{n\ge 0}X_{k_n}\quad\text{ and }\quad\widetilde X=\bigcup_{n\ge 0}\widetilde X_{n},
  $$ where $X_{n}=F_{n}\times C_{n+1}\times\cdots$
  and $\widetilde X_{n}=\widetilde F_{n}\times \widetilde C_{n+1}\times\cdots$.
Then the mappings
$$
 X_{k_n}\ni (f_{k_n}, c_{k_n+1},\dots)\mapsto (f_{k_n}, c_{k_n+1}\cdots c_{k_{n+1}},c_{k_{n+1}+1}\cdots c_{k_{n+2}},\dots)\in\widetilde X_n,
$$ 
$n\ge 0$, define a homeomorphism of $X$
onto $\widetilde X$.
This homeomorphism  intertwines $T$ with $\widetilde T$ and the $(C,F)$-equivalence relation on $X$ with the $(C,F)$-equivalence relations on $\widetilde X$.

\proclaim{Theorem 1.8} Let $T$ be the $(C,F)$-action of $G$ associated with a sequence $(C_n,F_{n-1})_{n\ge 1}$ satisfying \rom{(I)--(III)} and \thetag{1-3}.
Then $T$ is (measure theoretically) isomorphic to the  $(C,F)$-action $S$ of $G$ associated with a sequence $(C_n',F_{n-1}')_{n\ge 1}$ satisfying \rom{(I)--(III)} and \thetag{1-2}.
Moreover, the sequence $(F_n')_{n\ge 0}$ is a subsequence of $(F_n)_{n\ge 0}$.
In particular, if  $G=\Bbb Z$ and $T$ is of rank one then $S$ is also of rank one. 
\endproclaim

\demo{Proof}
Let $G=\{g_j\mid j\in\Bbb N\}$.
Denote by $X$ the $(C,F)$-space of  $T$.
Applying Proposition~1.4 and
passing, if necessary, to a $(C,F)$-subsequence we may assume without loss of generality
that for each $n$, there exists a subset $C_{n+1}'\subset C_{n+1}$ such that 
$$
\# g_jF_nC_{n+1}'\subset F_{n+1}, \ j=1,\dots, n,\quad\text{ and  }\quad\#C_{n+1}'/\# C_{n+1}>1-n^2.
$$
Then the sequence $(C_n',F_{n-1})_{n\ge 1}$ satisfies (I)--(III) and \thetag{1-2}.
We denote by $S$ the $(C,F)$-action of $G$ associated with $(C_n',F_{n-1})_{n\ge 1}$.
Let $Y$ stand for the corresponding $(C,F)$-space.
Then $Y=\bigcup_{n\ge 1}Y_n$, where $Y_n=F_n\times C_{n+1}'\times C_{n+2}'\times\cdots$.
By Borel-Cantelli lemma, for almost all $x\in X$, there is $n>0$ such that $x=(f_n,c_{n+1},c_{n+2},\dots)\in F_n\times C_{n+1}'\times C_{n+2}'\times \cdots\subset X_n$.
It is easy to see that the mapping 
$$
X\ni x\mapsto (f_n,c_{n+1},c_{n+2},\dots)\in Y_n\subset Y
$$
is a well defined (mod 0) isomorphism of $T$ with $S$.

The final claim of the theorem is obvious.
\qed
\enddemo

From Theorems 1.6, 1.8, Proposition~1.2  and Remark~1.3 we deduce the following corollary.

\proclaim{Corollary 1.9 (Minimal and almost minimal uniquely ergodic  models for rank-one actions)}
\roster
\item"\rom{(i)}"Every funny rank-one $\sigma$-finite measure preserving action of $G$ is measure theoretically isomorphic to  a strictly ergodic  topological $G$-action on a locally compact Cantor space.
\item"\rom{(ii)}"
Every funny rank-one $\sigma$-finite measure preserving action of $G$ is measure theoretically isomorphic to  an almost minimal uniquely ergodic\footnote{An almost minimal $G$-action is called uniquely ergodic if there is only one up to scaling non-atomic $\sigma$-finite $G$-invariant  Borel measure that is finite on every compact in the compliment to the fixed point of the action.}  topological $G$-action on a  compact Cantor space.
\endroster
\endproclaim

We note that the claims (i) and (ii) of the corollary are equivalent.

\remark{Remark \rom{1.10}} The authors of \cite{Dai--Si} (see also \cite{Bo--Wa}) introduced a concept of {\it normal} rank-one $\Bbb Z$-action.
This means that  (in the cutting-and-stacking construction of the action) at least one spacer is added above the highest  subtower
 for infinitely many cuts.
 As follows from Theorem~1.8,  every rank-one transformation is isomorphic to a normal one.
 We note however that passing to an isomorphic normal copy may ``destroy'' some other important properties of the rank-one  construction
such as the property of ``bounded cuts'' (which means that the sequence $(\# C_n)_{n=1}^\infty$ is bounded).
\endremark

We now state without proof (it is easy) a lemma which will be used in the next section.

\proclaim{Lemma 1.11}
Let $a,b\in F_n$.
Then $\mu(T_g[a]_n\cap[b]_n)>0$ if and only if 
$$
g\in \bigcup_{j>n}bC_{n+1}C_{n+2}\cdots C_jC_j^{-1}\cdots C_{n+2}^{-1}C_{n+1}^{-1}a^{-1}.
$$
Moreover,
$[a]_n\cap T_g^{-1}[b]_n=\sqcup[ad_{n+1}\cdots d_j]_{j}$, where the union is taken over all possible expansions of $g$ into the sum $g=bc_{n+1}c_{n+2}\cdots c_jd_j^{-1}\cdots d_{n+2}^{-1}d_{n+1}^{-1}a^{-1}$ with $c_l,d_l\in C_l$ for each $l$ and $c_j\ne d_j$.
\endproclaim

\remark{Remark \rom{1.12}}
The  $(C,F)$-construction of funny rank-one actions was introduced in \cite{dJ} and in \cite{Da1} in slightly different ways.
It was assumed in \cite{dJ} that the sequences  $(C_n,F_{n-1})_{n\ge 1}$  satisfy (I)--(III) plus an additional condition that $(F_n)_{n\ge 0}$ is a F{\o}lner sequence.
As we showed in Corollary~1.4, this additional condition is equivalent to~\thetag{1-3} in the case of finite measure preserving $(C,F)$-actions  (only such actions were under consideration in \cite{dJ}).
In \cite{Da1} we assumed that $(C_n,F_{n-1})_{n\ge 1}$  satisfy (I)--(III) and \thetag{1-2}.
An advantage of this approach is that  the associated $(C,F)$-actions are topological actions defined on locally compact spaces.
In this paper we introduced another condition \thetag{1-3} which is formally more general 
 than~\thetag{1-2}.
 However  it is shown in Theorem~1.8 that they define the very same class of the associated $(C,F)$-actions.
Thus the $(C,F)$-constructions from \cite{dJ} and \cite{Da1} are equivalent in the finite measure preserving case while the construction from \cite{Da1} is more general than that from \cite{dJ} in the i.m.p. case. 
In particular, for each non-amenable group $G$, there exist i.m.p. (and only i.m.p.) $(C,F)$-actions.
Thus we obtain a class of free ergodic conservative i.m.p. $G$-actions  whose orbit equivalence relations are  hyperfinite (see also \cite{Be--Go}, where another construction of such actions was given). 
\endremark

\head 2. Weak rational ergodicity and non-squashability  of funny rank-one Abelian actions
\endhead

Let $T$ be  an ergodic conservative measure preserving  action of an amenable group $G$ on a $\sigma$-finite measure space $(X,\goth B,\mu)$.
Fix an increasing  F{\o}lner sequence $(F_n)_{n\ge 0}$  in $G$.

\definition{Definition 2.1}
$T$ is called {\it  weakly rationally ergodic along $(F_n)_{n\ge 0}$} if there is a subset $Y$ of finite positive measure in $X$ such that
$$
\lim_{n\to\infty}\frac 1{a_n(Y)}\sum_{g\in F_n}\mu(A\cap T_gB)=\mu (A)\mu(B)\quad\text{for all }A,B\subset \goth B\cap Y,
$$
where $a_n(Y):=\sum_{g\in F_n}\frac{\mu(Y\cap T_g Y)}{\mu(Y)^2}$.
\enddefinition
We note that in the case where $G=\Bbb Z$ and $F_n=\{0,1,\dots,n-1\}$ we obtain the standard definition of the weak rational ergodicity \cite{Aa1}.
In the case where $G=\Bbb Z$ but $F_n=\{0,1,\dots,h_n-1\}$ for an increasing sequence
$(h_n)_{n=1}^\infty$, we obtain the notion  of the {\it subsequence weak rational ergodicity} \cite{Aa2}.

Our purpose in this section is to exhibit  a class of  funny rank-one actions of Abelian groups that are weakly rationally ergodic.

First of all we consider the case where $G=\Bbb Z$ and give a short proof of  the main result from the first version of \cite{Bo-Wa} that every rank-one transformation is subsequence weakly rationally ergodic\footnote{The two versions of  \cite{Bo-Wa} can be found in ArXiv.}  (see Theorem~2.4 below).
For that we need two auxiliary lemmata.

\proclaim{Lemma 2.2}
Let $T$ be the (rank-one) $(C,F)$-action of  $\Bbb Z$ associated with a sequence $(C_n,F_{n-1})_{n=1}^\infty$ with $F_n=\{0,\dots,h_n-1\}$ for each $n$ and  let $\mu$ be  the $(C,F)$-measure on the $(C,F)$-space $X$ such that $\mu([0]_0)=1$.
Then for all cylinders $A$ and $B$ in $X$,
$$
\lim_{n\to\infty}\frac{\sum_{k= 0}^{n-1}\mu(A\cap T_kB)}{a_n([0]_0)}=
\mu(A)\mu(B).
\tag2-1
$$
\endproclaim
\demo{Proof}
Fix $l\in\Bbb N$.
Let $A,B$ be $l$-cylinders, i.e. $A=[A']_l$ and $B=[B']_l$ for some $A',B'\subset F_l$.
Then 
$$
\mu(A\cap T_kB)=\sum_{a\in A',b\in B'}\mu([a]_l\cap T_k[b]_l)=\sum_{a\in A',b\in B'}\mu([0]_l\cap T_{k+b-a}[0_l]).\tag2-2
$$
In a similar way,
$$
\mu([0]_0\cap T_k[0]_0)=\sum_{c,c'\in C_1+\cdots +C_{l}}
\mu([0]_l\cap T_{k+c'-c}[0]_l).\tag2-3
$$
Since for each $m\in\Bbb Z$, 
the sequence $\sum_{k=1}^n\mu([0]_l\cap T_{k+m}[0]_l)$ of positive reals is equivalent\footnote{The equivalence of two sequences of reals means that the ratio of these sequences goes to $1$.} to the sequence $\sum_{k=1}^n\mu([0]_l\cap T_{k}[0]_l)$ as $n\to\infty$ , it follows from \thetag{2-2} and \thetag{2-3} that 
$$
\frac{\sum_{k=0}^{n-1}\mu(A\cap T_kB)}
{\sum_{k=0}^{n-1}\mu([0]_0\cap T_k[0]_0)}\to
\frac {\# A'\# B'}
{\#(C_1+\cdots+ C_l)^2}=
\mu([A']_l)\mu([B']_l),
$$
as desired.
\qed
\enddemo

\proclaim{Lemma 2.3}
Under the condition of the previous lemma,
for arbitrary subsets $A$ and $B$ in $[0]_0$, 
$$
\frac{\sum_{k=0}^{h_l-1}\mu(A\cap T_kB)}{\sum_{k=0}^{h_l-1}\mu([0]_0\cap T_k[0]_0)}\le
2\min(\mu(A),\mu(B)),
$$
where $h_l:=\# F_l$.
\endproclaim
\demo{Proof}
Let $C:=C_1+\cdots+C_l\subset \Bbb Z$.
For each subset $J\subset [0]_0$ and an element $c\in C$, we set $[c]_{l,J}:=[c]_l\cap J$.
Then $J=\bigsqcup_{c\in C}[c]_{l,J}$.
We now have 
$$
\align
\sum_{k=0}^{h_l-1}\mu(A\cap T_k B) &\le\sum_{k=0}^{h_l-1}\sum_{c,c'\in C}\mu([c]_{l,A}\cap T_k[c']_{l})\\
&=
\sum_{c,c'\in C}\sum_{k=0}^{h_l-1}\mu([c]_{l,A}\cap T_{c'} T_k[0]_{l})\\
&=
\sum_{c,c'\in C}\mu([c]_{l,A}\cap T_{c'} [0]_{0})\\
&\le
\sum_{c,c'\in C}\mu([c]_{l,A})\\
&=\mu(A)\#C.
\endalign
$$
In a similar way,
$$
\align
\sum_{k=0}^{h_l-1}\mu(A\cap T_k B) &\le
\sum_{k=0}^{h_l-1}\sum_{c,c'\in C}
\mu(T_{c+1-h_l}T_{-k}[h_l-1]_{l}\cap [c']_{l,B})\\
&=
\sum_{c,c'\in C}\mu( T_{c+1-h_l}     [0]_{0} \cap [c']_{l,B})\\
&\le
\sum_{c'\in C}\mu([c']_{l,B})\\
&=\mu(B)\#C.
\endalign
$$
Hence
$
\sum_{k=0}^{h_l-1}\mu(A\cap T_k B)\le \#C\cdot \min(\mu(A),\mu(B)).
$
On the other hand,
$$
\align
\sum_{k=0}^{h_l-1}\mu([0]_0\cap T_k [0]_0)\ge&
\sum_{c\ge c'\in C}\sum_{k=0}^{h_l-1}\mu([c]_{l}\cap T_k[c']_{l})\\
&=
\sum_{c\ge c'\in C}\sum_{k=0}^{h_l-1}\mu([c-c']_{l}\cap T_k[0]_{l})\\
&=
\sum_{c\ge c'\in C}\mu([c-c']_{l})\\
& =\frac{\# C+1}2,
\endalign
$$
and we are done.
\qed
\enddemo 

Approximating arbitrary pairs of Borel subsets in $[0]_0$ with cylinders, we deduce from Lemmata~2.2 and 2.3 the following assertion.

\proclaim{Theorem 2.4} 
For arbitrary Borel subsets $A$ and $B$ in $[0]_0$, 
$$
\lim_{\l\to\infty}\frac{\sum_{k=0}^{h_l-1}\mu(A\cap T_kB)}{\sum_{k=0}^{h_l-1}\mu([0]_0\cap T_k[0]_0)}=
\mu(A)\mu(B),
$$
i.e. $T$ is subsequence weakly rationally ergodic.
Hence
the rank-one transformations are non-squashable. 
\endproclaim

Now we pass to the general Abelian case.
Let $G$ be an  Abelian countable infinite discrete group $G$ and  let $T=(T_g)_{g\in G}$
be
 a $(C,F)$-action of $G$ associated with  $(C_n,F_{n-1})_{n\ge 1}$ 
 satisfying (I)--(III).
 In view of Theorem~1.8
 we may assume  without loss of generality that \thetag{1-2} holds.
 Moreover, we may assume without loss of generality that $0\in\bigcap_{n>0}F_n$.

We now prove an analogue of~Lemma~2.2.
Given a finite subset $D$ in $G$, we let $\nu_D:=\sum_{g\in D}\delta_g$, where $\delta_g$ is the Kronecker measure supported at $g\in G$.

\proclaim{Lemma 2.5}
The following are equivalent:
\roster
\item"\rom(i)"
For all cylinders $A,B\subset[0]_0$,
$$
\lim_{n\to\infty}\frac{\sum_{g\in F_{n}-F_{n}}\mu(A\cap T_gB)}{\sum_{g\in F_{n}-F_{n}}\mu([0]_0\cap T_g[0]_0)}=
\mu(A)\mu(B).\tag{2-4}
$$
\item"\rom(ii)"
For each $l>0$ and each $h\in F_l-F_l$, we have
$$
\sum_{g\in F_n-F_n+h}\mu([0]_l\cap T_g[0]_l)\sim_{n\to\infty}
\sum_{g\in F_n-F_n}\mu([0]_l\cap T_g[0]_l).
$$
\item"\rom(iii)"
For each $l>0$ and each $h\in F_l-F_l$, we have
$$
\lim_{n\to\infty}\lim_{L\to\infty}\frac{\nu_{C_{l+1}+\cdots+C_L}*\nu_{-C_{l+1}-\cdots-C_L}(F_n-F_n)}
{\nu_{C_{l+1}+\cdots+C_L}*\nu_{-C_{l+1}-\cdots-C_L}(F_n-F_n+h)}=1.
$$
\endroster
\endproclaim

\demo{Proof} (i)$\Rightarrow$(ii).
If (i)  holds then for each $l>0$ and  $f,v\in F_l$,
$$
1=\lim_{n\to\infty}\frac{\sum_{g\in F_{n}-F_{n}}\mu([f]_l\cap T_g[v]_l)}{\sum_{g\in F_{n}-F_{n}}\mu([0]_l\cap T_g[0]_l)}=
\lim_{n\to\infty}\frac{\sum_{g\in F_{n}-F_{n}}\mu([0]_l\cap T_{v+g-f}[0]_l)}{\sum_{g\in F_{n}-F_{n}}\mu([0]_l\cap T_g[0]_l)}
$$
and (ii) follows.

(ii)$\Rightarrow$(i) in the same way as in the proof of Lemma~2.2.

(ii)$\Leftrightarrow$(iii)
Fix a finite subset $S\subset G$.
Since \thetag{1-2} holds, for each $l>0$, there is $L>l$ such that
$$
S+C_{l+1}+\cdots+C_L\subset F_L.
$$
Therefore if $[c]_L\cap T_g[c']_L\ne\emptyset$ for some $c,c'\in C_{l+1}+\cdots+C_L$ and $g\in S$ then $c=g+c'$.
We now have
$$
\align
\sum_{g\in S}\mu([0]_l\cap T_g[0]_l)
&=\sum_{g\in S}\sum_{c,c'\in C_{l+1}+\cdots+ C_L}\mu([c]_L\cap T_g[c']_L)\\
&=\sum_{g\in S}\sum_{\{(c,c')\mid c,c'\in C_{l+1}+\cdots+ C_L, g=c-c'\}}\mu([c]_L\cap T_g[c']_L)\\
&=\sum_{g\in S}\sum_{\{(c,c')\mid c,c'\in C_{l+1}+\cdots+ C_L, g=c-c'\}}\frac 1{\# C_{1}\cdots\# C_L}\\
&=\frac {\nu_{C_{l+1}+\cdots+C_L}*\nu_{-C_{l+1}-\cdots-C_L}(S)}{\# C_{1}\cdots\# C_L}.
\endalign
$$
This implies the desired equivalence of (ii) and (iii).
\qed
\enddemo

Now we establish an analogue of Lemma~2.3.

\proclaim{Lemma 2.6} Suppose that there is $K>0$ such that $F_n-F_n\subset\bigsqcup_{j=1}^K d_{n,j}+F_n$ for some elements $d_{n,1},\dots, d_{n,K}\in G$.
Then
for arbitrary subsets $A$ and $B$ in $[0]_0$, 
$$
\frac{\sum_{g\in F_n-F_n}\mu(A\cap T_gB)}{\sum_{g\in F_n-F_n}\mu([0]_0\cap T_g[0]_0)}\le
K\min(\mu(A),\mu(B)).
$$
\endproclaim

\demo{Proof}
Let $C:=C_1+\cdots+C_n$.
For each subset $J\subset [0]_0$ and an element $c\in C$, we set $[c]_{n,J}:=[c]_n\cap J$.
Then $J=\bigsqcup_{c\in C}[c]_{n,J}$.
We now have 
$$
\align
\sum_{g\in F_n-F_n}\mu(A\cap T_g B) &\le\sum_{g\in F_n-F_n}\sum_{c,c'\in C}\mu([c]_{n,A}\cap T_g[c']_{n})\\
&\le
\sum_{c,c'\in C}\sum_{g\in F_n}\sum_{j=1}^K\mu([c]_{n,A}\cap T_{c'} T_{d_{n,j}+g}[0]_{n})\\
&=
\sum_{c,c'\in C}\sum_{j=1}^K\mu([c]_{n,A}\cap T_{c'} T_{d_{n,j}}[0]_{0})\\
&\le
K\sum_{c,c'\in C}\mu([c]_{n,A})\\
&=K\mu(A)\#C.
\endalign
$$
By symmetry, $\sum_{g\in F_n-F_n}\mu(A\cap T_g B)\le K\mu(B)\# C$.
On the other hand,
$$
\align
\sum_{g\in F_n-F_n}\mu([0]_0\cap T_g[0]_0)=&
\sum_{c, c'\in C}\sum_{g\in F_n-F_n}\mu([c]_{n}\cap T_g[c']_{n})\\
&\ge
\sum_{c,c'\in C}\sum_{g=c-c'}\mu([c]_{n}\cap T_g[c']_{n})\\
&=
\sum_{c, c'\in C}\mu([0]_{n})\\
& =\# C.
\endalign
$$
and we are done.
\qed
\enddemo

The corollary below follows from Lemmata~2.5 and 2.6 in the very same way as Theorem~2.4 follows from Lemmata~2.2 and 2.3.

\proclaim{Corollary 2.7} If the conditions of Lemmata~2.5  and 2.6 hold then \thetag{2-4} 
is satisfied for all Borel subsets $A,B\subset [0]_0$.
Hence $T$ is weakly rationally ergodic along the sequence $(F_n-F_n)_{n\ge 0}$.
In particular, $T$ is non-squashable. 
\endproclaim

We now give another sufficient condition for the weak rational ergodicity of Abelian $(C,F)$-actions.
Suppose that  for each $g\in G$, the   $(C,F)$-sequence $(C_n,F_{n-1})_{n\ge 1}$ satisfies  the following condition of  ``large holes''  in $C_{n+1}$:
$$
(g+F_n+F_n-F_n-F_n)\cap(C_{n+1}-C_{n+1})=\{0\}\tag{2-5}
$$
eventually.
Denote by $T=(T_g)_{g\in G}$
the associated $(C,F)$-action.
Without loss of generality we may assume that \thetag{1-2} holds.

\proclaim{Theorem 2.8}
If \thetag{2-5} holds then 
 \thetag{2-4} 
is satisfied for all Borel subsets $A,B\subset [0]_0$.
Thus $T$ is weakly rationally ergodic along $(F_n-F_n)_{n\ge 0}$.
\endproclaim

\demo{Proof}
We proceed in several steps.
 
 {\it Claim 1.} Given $g\in G$, we have  
 $$
 (g+F_n-F_n)\cap\sum_{j>0}(C_j-C_j)=\sum_{j=1}^n(C_j-C_j)\quad\text{ eventually.}
 $$
 Indeed, since $F_n+C_{n+1}\subset F_{n+1}$, it follows from \thetag{2-5} (with $n+1$ in place of $n$) that
  $$
 \bigg(g+\sum_{j=1}^2(F_n-F_n) +(C_{n+1}-C_{n+1})\bigg)\cap (C_{n+2}-C_{n+2})
  =\{0\}
 $$
 for all sufficiently large $n$.
From this and \thetag{2-5} we deduce that
 $$
 \bigg(g+\sum_{j=1}^2(F_n-F_n)\bigg)\cap \sum_{j=n+1}^{n+2}(C_j-C_j)=
 \bigg(g+\sum_{j=1}^2(F_n-F_n)\bigg)\cap (C_{n+1}-C_{n+1})=\{0\}
.
 $$
 By induction in $n$, we obtain that
$(g+\sum_{j=1}^2(F_n-F_n))\cap \sum_{j=n+1}^{\infty}(C_j-C_j)=\{0\}$ eventually.
Since $\sum_{j=1}^n(C_j-C_j)\subset F_n-F_n$,
we obtain that
$$
(g+F_n-F_n)\cap \sum_{j>0}(C_j-C_j)= (g+F_n-F_n)\cap \sum_{j=1}^n(C_j-C_j).
$$
 On the other hand, $-g+\sum_{j=1}^n C_j\subset F_n$ eventually in view of \thetag{1-2}.  Claim~1 follows.

{\it Claim 2.} For each $g\in G$, 
$$
\sum_{k\in g+F_n-F_n}\mu([0]_l\cap T_k[0]_l)=\sum_{k\in F_n-F_n}\mu([0]_l\cap T_k[0]_l)
\quad\text{ eventually.}
$$
Indeed, by Lemma~1.11, $\mu([0]_l\cap T_k[0]_l)>0$ if and only if $k\in\sum_{j>l}(C_j-C_j)$.
Claim~1 now yields that 
$$
(g+F_n-F_n)\cap \sum_{j>l}(C_j-C_j)=\sum_{j=l+1}^n(C_j-C_j)=(F_n-F_n)\cap \sum_{j>l}(C_j-C_j)
$$ 
and Claim~2 follows.

{\it Claim 3.} \thetag{2-4} holds for all cylinders $A,B$ in $X$.

\noindent This  follows from Lemma~2.5 and Claim~2.

{\it Claim 4.} For arbitrary subsets $A,B\subset[0]_0$,
$$
\frac{\sum_{g\in F_n-F_n}\mu(A\cap T_gB)}{\sum_{g\in F_n-F_n}\mu([0]_0\cap T_g[0]_0)}\le
\min\{\mu(A),\mu(B)\}.
$$
Indeed, let $C:=C_1+\cdots+C_n$ and $[c]_{n,A}:=[c]_n\cap A$ for each $c\in C$.
It follows from Lemma~1.11 and Claim~1 that
$$
\align
\sum_{g\in F_n-F_n}\mu(A\cap T_gB)& \le\sum_{g\in F_n-F_n}\sum_{c,c'\in C}\mu([c]_{n,A}\cap T_g[c']_{n})\\
&=\sum_{c,c'\in C}\sum_{g\in (F_n-F_n)\cap(c-c'+\sum_{j>n}(C_j-C_j))}\mu([c]_{n,A}\cap T_g[c']_{n})\\
&=\sum_{c,c'\in C}\mu([c]_{n,A}\cap T_{c-c'}[c']_{n})\\
&=\sum_{c,c'\in C}\mu([c]_{n,A})\\
&=\mu(A)\# C.
\endalign
$$
In a similar way, $\sum_{g\in F_n-F_n}\mu(A\cap T_gB) \le\mu(B)\# C$.
The same argument yields that $\sum_{g\in F_n-F_n}\mu([0]_0\cap T_g[0]_0)=\mu([0]_0)\# C=\# C$
and Claim~4 follows.

The assertion of the theorem follows now from Claims~3 and 4.
\qed
\enddemo

\head
3. Actions of finite funny rank  
\endhead

Let  $T$ be a measure preserving action of a discrete countable infinite group $G$
on a 
 standard $\sigma$-finite measure space $(X,\goth B,\mu)$.
 Fix $k>1$.
 
 \definition{Definition 3.1}
 We say that $T$ is of {\it funny rank at most $k$}
if there exist  $k$ sequences  $(B_n^{j})_{n\ge 0}$, $j=1,\dots,k$, of subsets of finite measure  in $X$ and  $k$ sequences  $(F_n^j)_{n\ge 0}$, $j=1,\dots,k$, of finite subsets in $G$ such that  
\roster
\item"(i)" 
for each $n\ge 0$ and $j\in\{1,\dots,k\}$, the subsets $T_gB_n^j$, $g\in F_n^j$, are pairwise disjoint and 
\item"(ii)" 
for each subset $B\in\goth B$ with $\mu(B)<\infty$,
$$
\lim_{n\to\infty}\inf_{F^j\subset F_n^j}\mu\bigg(B\triangle\bigsqcup_{j=1}^k\bigsqcup_{g\in F^j_n}T_gB_n^j\bigg)=0.
$$
\endroster
Without loss of generality we  may assume that $\lim_{n\to\infty}\mu(B_n^j)=0$ for each $j$.
The collection $\{T_gB_n^j\mid g\in F_n^j\}$ of subsets in $X$  is called the {\it $j$-th $T$-tower (of the $n$-th $T$-castle)}.
The subsets  $T_gB_b^j$, $g\in F_n^j$, are called {\it levels} of the
$j$-th $T$-tower. 
The {\it $n$-th $T$-castle} is the collection of $j$-th $T$-towers when $j$ runs the set $\{1,\dots,k\}$.
We say that the sequence of castles {\it refines} if for each $n>0$, every level of the $n$-th castle is a union of levels of the $(n+1)$-th castle.
The union of all levels of the $j$-th $T$-tower in the $n$-th $T$-castle  is denoted by $W_n^j$.
If $\mu(X\setminus\bigsqcup_{j=1}^kW_n^j)=0$ for all $n$,  then   we say that $T$ is {\it of funny rank at most $k$ without spacers}. 
If, moreover, there is $\delta>0$ such that $\mu(W_n^j)>\delta$ for all $n$ and $j$ then $T$ is of {\it exact funny rank  at most $k$}.
If, in addition,  there is $D>1$ such that $D^{-1}<\mu(B_n^i)/\mu(B_n^j)<D$ for all $n$ and $i,j\in\{1,\dots,k\}$ then we say that $T$ is {\it of balanced exact rank at most $k$}.
We say that $T$ is of {\it finite funny rank} if there is $k\ge 1$ such that $T$ is of funny rank  at most $k$.  
If $G=\Bbb Z$ and $F_n^j=\{0,1\dots,\# F_n^j-1\}$ then we obtain the standard definition of ``finite rank'', ``exact rank'', etc. (see \cite{Fe2} for details).
\enddefinition

The following lemma is standard.
We state it without proof.

\proclaim{Lemma 3.2}
Let $(B_n^j)_{n\ge 0}$ and $(F_n^j)_{n\ge 0}$ be as in Definition 3.1.
Given a sequence $(\epsilon_l)_{l=0}^\infty$ of positive reals tending to $0$, there is a sequence $(n_l)_{l=0}^\infty$ of positive integers  increasing to infinity and subsets $\widetilde B_l^j\subset B_{n_l}^j$ such that $\mu(B_{n_l}^j\setminus\widetilde B_l^j)<\epsilon_l$ for all $l$ and $j$,  the sequences $(\widetilde B_{l}^j)_{l\ge 0}$ and  $(F_{n_l}^j)_{l\ge 0}$, $j\in\{1,\dots,k\}$,  satisfy Definition~3.1 and the  sequence of the $n$-th $T$-castles corresponding to them refines. 
\endproclaim

Thus without loss of generality we may assume that if $T$
is of finite funny rank then the corresponding sequence of $T$-castles refines.

It is easy to see that every action of finite funny rank is conservative.
However  such actions  can be non-ergodic.
It is obvious that every action of finite funny rank without spacers is defined on a space with finite measure.

We now consider the case of $\Bbb Z$-actions of finite rank in more detail.

\proclaim{Lemma 3.3} Suppose that $T$ is an ergodic $\Bbb Z$-action of rank at most $k$.
If there is a Borel subset $A\subset X$  of finite strictly positive measure and $j\in\{1,\dots,k\}$  such that $\lim_{n\to\infty}\mu(A\cap W^j_n)=0$ then for each Borel subset $B\subset X$ of finite measure, $\lim_{n\to\infty}\mu(B\cap W_n^j)=0$. 
\endproclaim

\demo{Proof}
Let 
$$
\goth F:=\{C\in\goth B\mid \lim_{n\to\infty}\mu(C\cap W_n^j)=0\}.
$$
We claim that if $C\in\goth F$ and  $g\in \Bbb Z$ then $T_gC\in \goth F$. 
Indeed,
$$
\mu(T_gC\cap W_n^j)=\mu(C\cap T_{-g}W_n^j)\le
\mu(C\cap W_n^j) +\mu(W_n^j\triangle T_{-g} W_n^j)
$$
Since $\mu(W_n^j\triangle T_{-g} W_n^j)=2|g|\mu(B_n)\to 0$, it follows that
$T_gC\in\goth F$.

Now, given  a Borel subset $B\subset X$ of finite measure and $\epsilon>0$, there is $N>0$ such that
$\mu(B\setminus \bigcup_{m=0}^N T_mA)<\epsilon$.
Therefore
$$
\mu(B\cap W_n^j)\le \sum_{m=0}^N\mu( T_mA\cap W_n^j)+\epsilon
$$
Hence $B\in\goth F$.
\qed
\enddemo

We now show that  for the ergodic $\Bbb Z$-actions of finite rank, one can choose a sequence of  approximating $T$-castles (from Definition~3.1) that satisfy certain additional properties.

\proclaim{Proposition 3.4}
Let $T$ be an ergodic  $\Bbb Z$-action of rank at most $k_1$.
Then there is  $ k\le k_1$ and   sequences  $(B_n^j)_{n=1}^\infty$ and $(F_n^j)_{n=1}^\infty$, $1\le j\le k$ satisfying Definition~3.1 and such that the corresponding sequence
of $T$-castles refines and the following limits exist
$$
\lim_{n\to\infty}\mu\bigg(\bigg(\bigsqcup_{i=1}^k{B_0^i}\bigg)\cap W_n^j\bigg)=\delta_j>0,\quad j=1,\dots,k\tag3-1
$$
with $\sum_{j=1}^k\delta_j=\sum_{j=1}^k\mu(B_0^j)$.
\endproclaim

\demo{Proof}
We say that the sequence $(W_n^j)_{n=1}^\infty$ of $j$-th towers {\it vanishes} if 
$$
\lim_{n\to\infty}\mu\bigg(\bigg(\bigsqcup_{i=1}^{k_1}{B_0^i}\bigg)\cap W_n^j\bigg)=0.
$$
In follows from Lemma~3.3 that the vanishing towers do not really contribute into the property (ii) of Definition~3.1.
In other words, if we drop all the towers that vanish, the remaining sequences of towers still will satisfy  Definition~3.1.
Let $k$ be the number of the remaining sequences of towers.
Since they are non-vanishing,
there is a  subsequence of the associated  $T$-castles and strictly positive numbers $\delta_1,\dots,\delta_k$  such that \thetag{3-1} is satisfied.
Then, of course, 
$$
\sum_{j=1}^k\delta_j=\sum_{j=1}^k\mu(B_0^j).\tag3-2
 $$ 
 We need to modify the resulting sequence of $T$-castles to make it refining.
For that we apply Lemma~3.2 to find an increasing  sequence $(n_l)_{l\ge 0}$ of positive integers and subsets $(\widetilde B^j_l)_{l\ge 0}$ such that 
\roster
\item"(a)" $\widetilde B^j_l\subset B^j_{n_l}$ and
$\mu(B_{n_0}^j\setminus\widetilde B_0^j)<0.5k^{-1}\min_{1\le i\le k}\delta_i$ for all $l$ and $j$,
\item"(b)" 
 the sequences $(\widetilde B_{l}^j)_{l\ge 0}$ and  $(F_{n_l}^j)_{l\ge 0}$, $j\in\{1,\dots,k\}$  satisfy Definition~3.1 and
 \item"(c)"
  the  sequence of the $n$-th $T$-castles corresponding to them refines. 
\endroster
Passing to a further subsequence (of course, the subsequence satisfies the properties (a)--(c)) we may assume that there are reals $\widetilde\delta_1\ge 0,\dots,\widetilde\delta_k\ge 0$ such that
$\lim_{l\to\infty}\mu((\bigsqcup_{i=1}^k{\widetilde B_0^i})\cap \widetilde W_l^j)=\widetilde \delta_j$,
where $\widetilde W_l^j:=\bigsqcup_{g\in F^j_{n_l}}T_g \widetilde B_l^j$,
$j=1,\dots,k$.
It remains to show that the reals $\widetilde\delta_j$ are all strictly positive.
Of course, $\sum_{j=1}^k\widetilde\delta_j=\sum_{j=1}^k\mu(\widetilde B^j_0)$.
 If follows from (a) that $\widetilde\delta_j\le\delta_j$ for each $j$.
These inequalities, \thetag{3-2} and (a) yield that
 $$
 \sum_{j=1}^k\delta_j\ge \sum_{j=1}^k\widetilde\delta_j\ge \sum_{j=1}^k\delta_j-0.5\min_{1\le j\le k}\delta_j.
 $$
 This implies that $\widetilde\delta_j>0$ for each $j$, as desired.
 \qed
\enddemo

\proclaim{Corollary 3.5} Under the conditions of Proposition~3.4, the induced (finite measure preserving) transformation $(T_1)_{\bigsqcup_{i=1}^{k_1}B_0^i}$ is of exact rank at most $k$.
\endproclaim
\demo{Proof}
Suppose that the sequence of $T$-castles refines and satisfies \thetag{3-1}.
Then for each $n\in\Bbb Z_+$ and $j\in\{1,\dots,k\}$, the intersection of $W_n^j$  with  the space $\bigsqcup_{i=1}^{k}B_0^i$ is a tower of the induced transformation
$(T_1)_{\bigsqcup_{i=1}^{k}B_0^i}$.
The union of these towers, when $j$ runs $\{1,\dots,k\}$, is an $n$-th $(T_1)_{\bigsqcup_{i=1}^{k}B_0^i}$-castle.
The sequence of these castles refines and generates the entire $\sigma$-algebra of Borel subsets of $\bigsqcup_{i=1}^{k}B_0^i$.
Now \thetag{3-1} yields that $(T_1)_{\bigsqcup_{i=1}^{k}B_0^i}$ is exact. \qed
\enddemo

We recall that given a measure preserving transformation $S$ of a standard probability space $(Y,\nu)$, the {\it Koopman unitary operator} $U_S$ on  $L^2(Y,\nu)$ is defined by the formula $U_Sf:=f\circ S$. 

\proclaim{Lemma 3.6} Let  $S$ be an ergodic measure preserving transformation of a standard probability space $(Y,\nu)$.
Let $(A_n)_{n=1}^\infty$ be  a sequence of Borel subsets of $Y$ such that $\lim_{n\to\infty}\nu(A_n)=\delta>0$ and $\|U_S 1_{A_n}-1_{A_n}\|_2\to 0$.
Then for each Borel subset $B\subset Y$, $\nu(B\cap A_n)\to\nu(B)\delta$.
\endproclaim

\demo{Proof}
Since by von Neumann mean ergodic theorem, $\frac 1N\sum_{j=1}U_S^j1_B\to \mu(B)$ in the metric of $L^2(Y,\nu)$, we can find, for each $\epsilon>0$, a positive integer $N$ such that
$$
\epsilon\ge\bigg| \bigg\langle\frac 1N\sum_{j=1}^N U_S^j1_B,1_{A_n}\bigg\rangle-\nu(B)\nu(A_n)\bigg|
=\bigg| \bigg\langle1_B,\frac 1N\sum_{j=1}^N U_S^{-j}1_{A_n}\bigg\rangle-\nu(B)\nu(A_n)\bigg|
$$
for each $n\ge 0$.
Passing to the limit when $n\to\infty$ and using the condition of the lemma, we obtain that
$\epsilon\ge\limsup_{n\to\infty}|\nu(B\cap A_n)-\nu(B)\delta|$.
Hence $\nu(B\cap A_n)\to\nu(B)\delta$.
\qed
\enddemo

Applying Lemma~3.6 we  refine Proposition~3.4 in the following way.

\proclaim{Corollary 3.7}
Passing to a further subsequence in $(B_n^j, F_n^j)_{n=1}^\infty$ from the statement of Proposition~3.4 and normalizing $\mu$ such that $\mu((\bigsqcup_{l=1}^k{B_0^l})=1$ we may assume without loss of generality that the following holds: \thetag{3-1} and
$$
\lim_{n\to\infty}\mu\bigg(A\cap W_n^i\cap W_{n+1}^j\bigg)=\mu(A)\delta_i\delta_j,\quad i,j=1,\dots,k,
$$
for each Borel subset $A\subset\bigsqcup_{l=1}^k{B_0^l}$.
Moreover, there exist limits
$$
\Lambda_i:=\lim_{n\to\infty}\frac{\mu(B_n^i)}{\sum_{l=1}^k\mu(B_n^l)}
$$
and for each level $D_n^i$ of the $i$-th tower of the $n$-th $T$-castle such that $D_n^i\subset\bigsqcup_{l=1}^k{B_0^l}$,
$$
\lim_{n\to\infty}\frac{\mu( D_n^i\cap W_{n+1}^j)}
	{\mu(( \bigsqcup_{l=1}^k D_n^l)\cap W_{n+1}^j)}
	=\Lambda_i,\qquad i=1,\dots,k.
	$$
	Of course, $\sum_{i=1}^k\Lambda_i=1$.
\endproclaim

\demo{Proof}
It suffices to note  
 that 
 \roster
 \item"---"
 the reals $\delta_1,\dots,\delta_{k_1}$ do not depend on passing to a subsequence of approximating castles,
 \item"---"
 the intersection of $W_n^j$ with the set $\bigsqcup_{l=1}^k{B_0^l}$ is a tower of the ergodic probability preserving induced transformation $(T_1)_{\bigsqcup_{l=1}^k{B_0^l}}$
 \endroster
 and apply Lemma~3.6.
 \qed
\enddemo

\head 4. $(C,F)$-construction  of actions of  finite funny rank
\endhead

We now give a constructive definition for actions of $G$ of finite funny rank.

We set $G':=G\times\{1,\dots,k\}$ and $G'':=\{1,\dots,k\}\times G\times\{1,\dots,k\}$.
Given $\boldsymbol c\in  G''$, we denote by $ c$ the  element in $G$ such that $\boldsymbol c=(i, c,j)$ for some $i,j\in\{1,\dots,k\}$. 
In a similar way, given $\boldsymbol a\in G'$, we denote by $ a$ the element in $G$ such that
$\boldsymbol a=( a,i)$ for some $i\in\{1,\dots,k\}$.
We will consider  $G'$ as a left $G$-space, where   $G$ acts by the formula $g\cdot ( f,i):=(gf,i)$.
Given a subset $A$ of $G'$ and $i\in\{1,\dots,k\}$, we  let $A^i:=A\cap (G\times\{i\})$.
In a similar way, given a subset $C$ of $G''$ and  $i,j\in\{1,\dots,k\}$, we denote by $C^{i,j}$ the intersection 
$
C\cap(\{i\}\times G\times\{j\}).
$
Let $(f,i)\in G'$ and $(k,g,l)\in G'$.
If   $i=k$ we define a ``product'' $(f,i)*(k,g,l)$  by setting 
$$
(f,i)*(k,g,l):=(fg,l)\in G'.
$$ 
For arbitrary subsets $A\subset G'$ and $C\subset G''$, we let $A*C$ be the set of all products $\boldsymbol a*\boldsymbol c$, where $\boldsymbol a\in A$, $\boldsymbol c\in C$ and $\boldsymbol a*\boldsymbol c$ is defined.
We reduce the notation $A*\{\boldsymbol c\}$ to $A* \boldsymbol c$.
In a similar way one can define a product of two elements of $G''$:
$$
(i,c,j)*(k,c',k):=(i,cc',k)\in G''\quad \text{if}\quad j=k.
$$
Hence  the product $C*C'$ of two subsets  $C,C'$ in $G''$ is also well defined.

Suppose we are given two sequences of finite subsets $(F_n)_{n\ge 0}$ in $G'$ and $(C_n)_{n\ge 1}$ in $ G''$ such that the following conditions hold for each $n$:

\roster
\item"(I)"
$F_0=\{1\}\times\{1,\dots,k\}$, $\sum_{j=1}^k\# C_{n}^{i,j}>1$  and $(1,i)\in F_n$ for each $i$,
\item"(II)"
$F_n *C_{n+1}\subset F_{n+1}$, 
\item"(III)"
$F_n *\boldsymbol c\cap F_n*\boldsymbol c'=\emptyset$ if $\boldsymbol c\ne \boldsymbol c'\in C_{n+1}$.
\endroster
We let for each $n\ge 0$,
 $$
 \align
 X_n :=   \{(\boldsymbol f_n,\boldsymbol c_{n+1},\boldsymbol c_{n+2}\dots) &\in F_n\times C_{n+1}\times C_{n+2}\times\cdots \\
&\mid \boldsymbol f_n*\boldsymbol c_{n+1}*\cdots*\boldsymbol c_{l}\text{ is well defined for each $l>n$}\}.
 \endalign
 $$
Then $X_n$ is a perfect subset of the compact Cantor  space $F_n\times C_{n+1}\times C_{n+2}\times\cdots $ (endowed with Tikhonov's topology).
Hence $X_n$ is itself a compact Cantor space.
The map
$$
X_n\ni(\boldsymbol f_n,\boldsymbol c_{n+1},\boldsymbol c_{n+2}\dots) \in (\boldsymbol f_n*\boldsymbol c_{n+1},\boldsymbol c_{n+2}\dots) \in X_{n+1}
$$
is  a topological embedding of $X_n$ into $X_{n+1}$.
Therefore we will consider $X_n$ as a (clopen) subset of $X_{n+1}$.
We now define $X$ to be the topological inductive limit of the increasing sequence $X_1\subset\ X_2\subset\cdots$ of compact Cantor spaces.
Then $X$ is a locally compact Cantor space\footnote{$X$ is compact if and only if there is $N>0$ with $F_{n+1}=F_n*C_{n+1}$ for all $n>N$.}.
We call it the {\it $(C,F)$-space associated with the sequence} $(C_n,F_{n-1})_{n\ge 1}$.
Given a subset $A$ of $F_n$, we let 
$$
[A]_n:=\{(\boldsymbol f_n,\boldsymbol c_{n+1},\dots)\in X_n\mid \boldsymbol f_n\in A\}\subset X.
$$
and call this set an {\it $n$-cylinder}.
It is a compact open subset of $X$.
Conversely, every compact open subset of $X$ is a cylinder. 
The set of all cylinders is a base of the topology in $X$.
It is easy to see that
$$
[A]_n\cap[B]_n=[A\cap B]_n,\quad 
[A]_n\cup[B]_n=[A\cup B]_n\quad\text{and}\quad [A]_n=[A*C_{n+1}]_{n+1}
$$ 
for all $A,B\subset F_n$ and  $n\ge 0$.
For brevity we will write $[\boldsymbol f]_n$ for $[\{\boldsymbol f\}]_n$, $\boldsymbol f\in F_n$.

Let $\Cal R$ denote the {\it tail equivalence relation} on $X$.
 This means that the restriction of $\Cal R$ to $X_n$ is the tail equivalence relation on $X_n$ for each $n \ge 0$.
If the following condition is satisfied (in addition to (I)--(III)):
 \roster
 \item"(IV)" $\# C_n^{i,j}>1$ for all $i,j\in\{1,\dots,k\}$ and $n>0$
 \endroster
then
 $\Cal R$ is {\it minimal}, i.e. the $\Cal R$-class of every point in $X$ is dense in $X$.
 
We now investigate the problem of existence and uniqueness of $\Cal R$-invariant Radon measures   on $X$.
Of course, such a measure $\mu$ is completely determined by its values on the cylinders.
In turn, every cylinder is a disjoint union of some ``elementary'' cylinders $[\boldsymbol f]_n$, where $\boldsymbol f$ runs $F_n$.
Since the elementary cylinders $[\boldsymbol f]_n$ and $[\boldsymbol g]_n$ are of the same measure  whenever $\boldsymbol f,\boldsymbol g\in F_n^i$ for some $i\in\{1,\dots,k\}$ (this fact is equivalent to the $\Cal R$-invariance of $\mu$), we obtain that $\mu$ is determined uniquely by its values  
$$
\lambda_n^i:=\mu([(1,i)]_n),\quad i=1,\dots, k, \quad\text{for all } n\ge 0.
\tag4-1
$$
Thus, we obtain a sequence of  vectors $\lambda_n:=\pmatrix \lambda^1_n\\ \vdots\\ \lambda^k_n\endpmatrix\in\Bbb R^k_+$, $n\ge 0$.
There is a consistency condition for these vectors.
Indeed, the  property 
$$
\mu([\boldsymbol f]_n)=\sum_{\boldsymbol c\in C_{n+1}}\mu([\boldsymbol f*\boldsymbol c]_{n+1}),\quad  \boldsymbol f\in F_n,
$$ 
can be rewritten as
$$
\lambda_n=r_{n+1}\lambda_{n+1},\tag 4-2
$$
where
 $r_n=(r_n^{i,j})_{1\le i,j\le k}$  is a  $k\times k$ integer matrix defined by setting $r_n^{i,j}:=\#(C^{i,j}_n)$.
Conversely, given a sequence $(\lambda_n)_{n\ge 0}$ of positive vectors in $\Bbb R^k$ satisfying \thetag{4-2}, there is a unique $\Cal R$-invariant  Radon measure $\mu$ on $X$ satisfying \thetag{4-1}.
Indeed, we define $\mu$ by setting
 for each $n\ge 0$ and a subset $A\subset F_n$,
$$
\mu([A]_n):=
\sum_{i=1}^k\#(A^i)\lambda_n^i
.\tag4-3
$$
We call $\mu$  the {\it $(C,F)$-measures on $X$ associated with the sequence} $(\lambda_n)_{n\ge 0}$.
It is easy to see that $\mu$ is finite if and only if 
$$
\lim_{n\to\infty}\sum_{j=1}^k\lambda_{n}^i\#(F_{n}^i)<\infty.\tag4-4
$$

Thus we obtain the following proposition.

\proclaim{Proposition 4.2} The formula \thetag{4-3} establishes  an affine  one-to-one correspondence between the set of $\Cal R$-invariant Radon measures  on $X$ and  the projective limit of the  sequence 
$$
\Bbb R_+^k@<{r_1}<<\Bbb R_+^k@<{r_2}<<\Bbb R_+^k@<{r_3}<<\cdots.
$$
\endproclaim

Unlike the funny rank-one case considered in Section~1, in the general case there can be several mutually disjoint ergodic $\Cal R$-invariant Radon measures on $X$.
Now we provide a simple sufficient condition on the sequence $(r_n)_{n=1}^\infty$ under which  $\Cal R$  is {\it uniquely ergodic}, i.e. there exists a unique  (up to scaling) ergodic $\Cal R$-invariant non-trivial Radon measure on $X$.

\proclaim{Proposition 4.3} 
Suppose that  there exist nonnegative reals $\Lambda_1,\dots,\Lambda_k$ such that for an increasing subsequence  of integers $n_p\to\infty$, there exist limits
$$
\lim_{p\to\infty} \frac{r^{i,l}_{n_p}}{\sum_{j=1}^kr^{j,l}_{n_p}}=
\Lambda_i \qquad\text{for all $i,l\in\{1,\dots,k\}$}. \tag4-5
$$
Then  $\Cal R$
is uniquely ergodic.
 \endproclaim

\demo{Proof}
In view of Proposition~4.2,
it suffices to show that the intersections of the cones $r_{n_p}(\Bbb R^k_+)$ with the simplex $\Delta:=\{(z_1,\dots,z_k)\in\Bbb R^k_+\mid z_1+\dots+z_k=1\}$ shrink to a singe point as $p\to\infty$.
The cone $r_{n_p}(\Bbb R^k_+)$
is  generated by $k$ rays passing through the following vectors $(r_{n_p}^{1,1},\dots,r_{n_p}^{k,1})$,\dots,
$(r_{n_p}^{1,k},\dots,r_{n_p}^{k,k})\in\Bbb R^k_+$.
Therefore~\thetag{4-5} yields that $ r_{n_p}(\Bbb R^k_+)\cap\Delta\to \{(\Lambda_1,\dots,\Lambda_k)\}$ (in Hausdorff metric) as $p\to\infty$.
\qed
\enddemo

\remark{Remark \rom{4.4}} We also note that $X_0$ intersects each $\Cal R$-orbit infinitely many 
times.
The map $\mu\mapsto\mu\restriction X_0$ is an affine isomorphism of the simplex  of Radon $\Cal R$-invariant measures on $X$  which equal $1$ on $X_0$ onto the simplex of $(\Cal R\restriction X_0)$-invariant probability  Borel measures on $X_0$.
Indeed, every $(\Cal R\restriction X_0)$-invariant measure extends uniquely to an $\Cal R$-invariant measure on $X$ via \thetag{4-3} and \thetag{4-1}.
Therefore this isomorphism maps the ergodic  $\Cal R$-invariant Radon measures which equal $1$ on $X_0$ onto  the ergodic  $(\Cal R\restriction X_0)$-invariant probability Borel measures.
In particular, $\Cal R$ is uniquely ergodic if and only if so is $\Cal R\restriction X_0$.
\endremark

We now define an action of $G$ on $X$ (or, more rigorously,  on a subset of $X$).
Given $g\in G$, let
$$
X_n^g:=\{(\boldsymbol f_n,\boldsymbol c_{n+1},\boldsymbol c_{n+2}\dots)\in X_n\mid g\cdot \boldsymbol f_n\in F_n\}.
$$
Then $X_n^g$ is a compact open subset of $X_n$ and $X_n^g\subset X_{n+1}^g$.
Hence the union $X^g:=\bigcup_{n\ge 0}X_n^g$ is  a well defined  open subset of $X$.
Let $X^G:=\bigcap_{g\in G}X^g$.
Then $X^G$ is a $G_\delta$-subset of $X$.
It is Polish in the induced topology.
Given $x\in X^G$ and $g\in G$, there is $n>0$ such that $x=(\boldsymbol f_n,\boldsymbol c_{n+1},\dots)\in X_n$ and $g\cdot \boldsymbol f_n\in F_n$.
We now let 
$$
T_gx=(g\cdot \boldsymbol f_n, \boldsymbol c_{n+1},\dots)\in X_n\subset X.
$$
It is standard to verify that $T_gx\in X^G$, the map $T_g:X^G\ni x\mapsto T_gx\in X^G$ is a homeomorphism of $X^G$ and $T_gT_{g'}=T_{gg'}$ for all $g,g'\in G$.
Thus $T:=(T_g)_{g\in G}$ is a continuous action of $G$ on $X^G$.

\definition{Definition 4.5}
We call $T$ {\it the $(C,F)$-action of $G$ associated with the sequence $(C_n,F_{n-1})_{n\ge 0}$}.
\enddefinition
The $(C,F)$-actins are free.
It is obvious that  $X^G$ is $\Cal R$-invariant and the $T$-orbit equivalence relation is the restriction of $\Cal R$ to $X^G$.
Hence for each ergodic $(C,F)$-measure $\mu$, either $\mu(X^G)=0$ or  $\mu(X\setminus X^G)=0$ and $T$ preserves  $\mu$.

We now state an analogue of Proposition~1.2.

\proclaim{Proposition 4.6}
 $X^G=X$ if and only if for each $g\in G$ and $n>0$, there is $m>n$ such that
$$
g\cdot F_n*C_{n+1}*C_{n+2}*\cdots * C_{m}\subset F_m.\tag 4-6
$$
\endproclaim

Thus, in this case the $(C,F)$-action is defined on the entire (locally compact) space $X$.
We do not give the proof of this lemma because it is an obvious slight modification of the proof of Proposition~1.2.
As in the rank-one case considered in Section~2, if $X$ is not compact then $T$ extends continuously to the one-point compactification of $X$.
If (IV) is satisfied then the extended action is almost minimal.

The following lemma is a counterpart of Proposition~1.4.

\proclaim{Proposition 4.7}
Let $\mu$ be a $(C,F)$-measure on $X$ associated with a sequence $(\lambda_n)_{n\ge 0}$.
\roster
\item"\rom{(i)}"
 $\mu(X\setminus X^G)=0$
 if and only if for each $g\in G$ and $n>0$,
$$
\lim_{m\to\infty}{\sum_{i=1}^k\#(A_{m,n}^i)\lambda_m^i}={\sum_{i=1}^k\#( F_n^i)\lambda_n^i},
\tag4-7
$$
where $A_{m,n}:=(g\cdot F_n*C_{n+1}*\cdots *C_{m})\cap F_m\subset G'$.
\item"\rom{(ii)}" 
If $\mu(X)<\infty$ and there exist limits $\gamma_i:=\lim_{n\to\infty}\frac{\lambda_n^i\# F_n^i}{\sum_{j=1}^k \lambda_n^j\# F_n^j}$, $i=1,\dots,k$\footnote{Passing to a $(C,F)$-subsequence, we can assume without loss of generality that the latter condition holds always.}, 
then  $\mu(X\setminus X^G)=0$ if and only if the sequence $(F_n^i)_{n=1}^\infty$ is F{\o}lner in $G$ for each $i$ such that $\gamma_i\ne 0$.
In particular, $G$ is amenable.
\endroster
\endproclaim

\demo{Proof}
We only prove (ii).
We note that  $\mu(X\setminus X^G)=0$ if and only if for each $g\in  G$, $\mu(X_n^g)/\mu(X_n)\to 1$. 
Since
$$
\frac{\mu(X_n^g)}{\mu(X_n)}=\frac{\sum_{i=1}^k\#(g\cdot F_n\cap F_n)^i\lambda_n^i}{\sum_{i=1}^k\# F_n^i\lambda_n^i}=\sum_{i=1}^k\frac{\#(g\cdot F_n\cap F_n)^i}{\# F_n^i}
\cdot\frac{\lambda_n^i\# F_n^i}{\sum_{j=1}^k\lambda_n^j\# F_n^j},
$$
it follows from this and the condition of the proposition that
$$
\sum_{i=1}^k\frac{\#((g\cdot F_n)^i\cap F_n^i)}{\# F_n^i}
\gamma_i\to 1.\tag4-8
$$
Since $\gamma_i\ge 0$ for each $i$ and $\sum_{i=1}^k{\gamma_i}=1$, it follows that \thetag{4-8} is equivalent to the following claim:
$$
\frac{\#((g\cdot F_n)^i\cap F_n^i)}{\# F_n^i}\to 1\qquad\text{for each $i$ such that $\gamma_i\ne 0$.}
$$
\qed
\enddemo

We note that if $G$ is amenable, $\mu(X)<\infty$ and $(F_n^i)_{n\ge 0}$ is a F{\o}lner sequence in $G$ for each $i\in\{1,\dots,k\}$ then~\thetag{4-7} holds.

From now on we assume that \thetag{4-7} is satisfied.
Then $(X,\mu,T)$ is a (well defined) measure preserving dynamical system.
It is easy to see that it is conservative.

By  analogy with the case of $(C,F)$-actions of rank one, we may assume without loss of generality that the sequences $(F_n)_{n\ge 0}$ and $(C_n)_{n\ge 1}$ satisfy the following condition
\roster
\item"(V)" 
$\{(1,i)\mid 1\le i\le k\}\subset\bigcap_{n=0}^\infty F_n$ and $\{(i,1,i)\mid 1\le i\le k\}\subset\bigcap_{n=1}^\infty C_n$
\endroster
in addition to (I)--(III).

We claim that each $(C,F)$-action (defined in this section)
 is of funny rank at most $k$.
Indeed, the sequences $([(1,1)]_n)_{n\ge 0}, \dots, ([(1,k)]_n)_{n\ge 0}$ and $(F_n^1)_{n\ge 0},\dots,
(F_n^k)_{n\ge 0}$ satisfy Definition~3.1.
We also note that $T_g[\boldsymbol f]_n=[g\cdot \boldsymbol f]_n$  (up to a $\mu$-null subset) whenever $\boldsymbol f,g\cdot \boldsymbol f\in F_n$.

We collect   some of the above  results on  $(C,F)$-actions  of finite funny rank in the following theorem.

\proclaim{Theorem 4.8}
Given a sequence $(C_n,F_{n-1})_{n\ge 1}$ satisfying  (I)--(V) and \thetag{4-5}, there is a locally compact Cantor  space $X$ and a countable equivalence relation $\Cal R$ on $X$ such that
\roster
\item"(i)"
 every $\Cal R$-class is dense in $X$,
\item"(ii)"
there is only one (up to scaling) $\Cal R$-invariant non-trivial $\sigma$-finite Radon measure $\mu$ on $X$,
\item"(iii)"
$\mu$ is finite if and only if \thetag{4-4} and \thetag{4-1} are satisfied,
\item"(iv)" 
there is a free topological  $G$-action $T$ on an $\Cal R$-invariant $G_\delta$-subset $X^G$ of $X$ such that the $T$-orbit equivalence relation is the restriction of $\Cal R$ to $X^G$,
\item"(v)"
$X^G=X$ if and only if \thetag{4-6} is satisfied,
\item"(vi)"
$\mu(X\setminus X^G)=0$ if and only if \thetag{4-7} is satisfied.
If  \thetag{4-7} is not satisfied then $\mu(X^G)=0$.
\item"(vii)"
If  $\mu(X)<\infty$ and $(F_n^j)_{n\ge 0}$ is a F{\o}lner sequence in $G$ for each $j=1,\dots,k$ then \thetag{4-7} is satisfied.
\item"(viii)"
Under  \thetag{4-7}, the dynamical system $(X,\mu,T)$ is ergodic, conservative and of funny rank at most $k$.
\endroster
\endproclaim

The converse to Theorem~4.8(viii) also holds.

\proclaim{Theorem 4.9 ($(C,F)$-models for finite funny rank  actions)} If $T$ is a $G$-action of funny rank at most $k$ on a standard $\sigma$-finite measure space $(X,\goth B,\mu)$ then $T$ is isomorphic to a $(C,F)$-action of $G$ and the corresponding isomorphism maps $\mu$ to a $(C,F)$-measure.
\endproclaim

We do not give a proof of this theorem because it is an obvious   modification of the proof of Theorem~1.6.

Let $T$ be the  $(C,F)$-action associated with  a sequence $(C_n,F_{n-1})_{n>0}$ satisfying (I)--(III) and \thetag{4-7}.
Given an increasing  sequence $(k_n)_{n\ge0}$ of non-negative integers with $k_0=0$, we let
$\widetilde F_n:=F_{k_n}$ and $\widetilde C_{n}:=C_{k_{n-1}+1}*C_{k_{n-1}+2}*\cdots * C_{k_{n}}$.
Then  the sequence $(\widetilde C_n,\widetilde F_{n-1})_{n>1}$  also satisfies (I)--(III) and \thetag{4-7}.
By an analogy with the rank-one $(C,F)$-actions, we call $(\widetilde C_n,\widetilde F_{n-1})_{n>1}$ a {\it $(C,F)$-subsequence} of $( C_n, F_{n-1})_{n>1}$.
 The $(C,F)$-action associated with it is
 canonically isomorphic to $T$.

We also state a higher rank analogue of Theorem~1.8.

\proclaim{Theorem 4.10} Let $T$ be the $(C,F)$-action of $G$ associated with a sequence $(C_n,F_{n-1})_{n\ge 1}$ satisfying \rom{(I)--(III)} and \thetag{4-7}.
Then $T$ is (measure theoretically) isomorphic to the  $(C,F)$-action $S$ of $G$ associated with a sequence $(C_n',F_{n-1}')_{n\ge 1}$ satisfying \rom{(I)--(III)} and \thetag{4-6}.
Moreover, the sequence $(F_n')_{n\ge 0}$ is a subsequence of $(F_n)_{n\ge 0}$.
In particular, if  $G=\Bbb Z$ and $T$ is of rank at most $k$ then $S$ is also of rank at most $k$. 
\endproclaim

We omit the proof of this theorem because it is a slight modification of the proof of Theorem~1.8.

\remark{Remark  \rom{4.11}}
In the case where $G=\Bbb Z$, $F_n^i=\{0,\dots,h_n^i-1\}$ and the $(C,F)$-space $X$ is compact (i.e. when the associated $(C,F)$-action $T$ has rank at most $k$ without spacers and the corresponding sequence of $T$-castles refines) we can associate an ordered Bratteli diagram with $T$.
We refer to \cite{Du} and \cite{Be--So} for the definitions related to Bratteli diagrams and Bratteli-Vershik systems.
We define a graded vertex set $V=(V_n)_{n\ge -1}$ and a graded edge set $E:=(E_n)_{n\ge 0}$ in the following way.
We let  $V_{-1}:=\{0\}$,  $V_n:=\{1,\dots,k\}$  for $n\ge 0$ and $E_0:=F_0$ and $E_n:=C_n$ for $n\ge 1$.
More precisely, we consider $C_n^{i,j}$ as the set of edges connecting the vertex $i\in V_{n-1}$ with the vertex $j\in V_n$, $n\ge 1$. 
The corresponding  graded graph $(V,E)$ is called a {\it Bratteli diagram}.
We now define an order relation on it.
For that we  define for each $n\ge 1$ and $j\in V_n$, a linear order $\succ$ on the set $\bigsqcup_{i=1}^kC_n^{i,j}$ by setting:
 $\boldsymbol c\succ \boldsymbol d$ if $ c\ge d$.
 Thus, we obtain an ordered Bratteli diagram $(V,E,\succ)$.
It is straightforward to verify that the map $\phi: X_0\ni(f_0,c_1,c_2,\dots)\mapsto(f_0,c_1,c_2,\dots)$ identifies $X_0$ with the Bratteli compactum   $Y$ associated with $(V,E)$.
Moreover, $\phi$ intertwines $T$ with the Bratteli-Vershik  $\Bbb Z$-action associated with $(V,E,\succ)$ and $\phi\times\phi$ maps bijectively the $(C,F)$-equivalence relation on $X$ onto the tail equivalence relation on $Y$. 
\endremark

\head 5. 
Weak rational ergodicity of transformations  of balanced finite  rank 
\endhead

Let $T$ be an ergodic  i.m.p.  $\Bbb Z$-action of finite rank.
Let $(X,\mu)$ stand for  the space of this action.
By Theorem~4.9, $T$ is isomorphic to a $(C,F)$-action of $\Bbb Z$ associated with sequences
$(C_n)_{n>0}$ and $(F_n)_{n\ge 0}$ satisfying (I)--(III) and (V) from Section~4.
Thus we consider $X$ as a $(C,F)$-space and $\mu$ as a $(C,F)$-measure associated with a sequence $(\lambda_n)_{n=1}^\infty$.
Since $T$ is ergodic,  we may apply Proposition~3.4 and Corollary~3.7.
``Translating'' their assertions into the language of $(C,F)$-systems,
 we may assume without loss of generality that  
\roster
\item"$(\alpha)$"
$\mu(X_0)=1$ and
the induced transformation $(T_1)_{ X_0}$ is an ergodic probability preserving transformation of exact rank at most $k$,
\item"$(\beta)$"
 there exist  limits    $\lim_{n\to\infty}\mu(X_0\cap[F_n^j]_n)=\delta_j>0$ for each $j\in\{1,\dots,k\}$, and $\sum_{j=1}^k\delta_j=1$,
 \item"$(\gamma)$"
 there exist limits 
 $$
\lim_{n\to\infty}\frac{\lambda^i_n}{\sum_{j=1}^k\lambda^j_n}
=\lim_{n\to\infty}\frac{r_{n}^{i,j}}{\sum_{l=1}^k r_{n}^{l,j}}=\Lambda_i\ge 0
$$ 
for all $i,j\in\{1,\dots,k\}$ and $\sum_{i=1}^k\Lambda_i=1$.
\endroster

We note that id $\Lambda_i>0$ for all $i$ then (IV) from Section~4 is satisfied.
Hence $T$ is a minimal $(C,F)$-action.

It is well known that each transformation of rank one is isomorphic to a {\it tower transformation}\footnote{We recall a standard definition. Given a measure preserving transformation $S$ on a standard measure space $(Y,\nu)$ and a Borel map $f:Y\to\Bbb Z_+$, we define a new dynamical system $(X,\mu,T)$ by setting $X=\{(y,i)\in Y\times\Bbb Z\mid 0\le i<f(y)\}$, $d\mu(y,i)=d\mu(x)$ and $T(y,i)=(y,i+1)$ if $i+1<f(y)$ and $T(y,i)=(Sy,i)$ if $i+1=f(y)$.
Then $T$ is called the transformation built under $f$ over $S$.} built under certain ``spacer map'' over  an ergodic transformation with pure point spectrum.
We now extend this result to the transformations of finite rank.
For that we first introduce  a sequence $(s_n)_{n\ge 0}$ of auxiliary mappings from $X_0$ to $\Bbb Z_+$ by setting
$$
s_n(\boldsymbol f_0,\boldsymbol c_1,\boldsymbol c_2,\dots):=\max\{t\ge 0\mid t\cdot \boldsymbol f_0*\boldsymbol c_1*\cdots *\boldsymbol c_{n}\in F_{n}\setminus (F_{0}*C_1*\cdots  *C_{n})\},
$$
where $\boldsymbol f_0\in F_0$ and $\boldsymbol c_j\in C_j$, $n\ge 0$.
Of course, $s_n$ is  continuous for each $n\ge 0$ and $s_0\le s_1\le\cdots$.
We let 
$$
C_n^{\max}:=\bigg\{ \boldsymbol c\in  \bigsqcup_{i=1}^kC_{n}^{i,j}\, \bigg| \, c=\max_{ \boldsymbol d\in \bigsqcup_{i=1}^kC_{n}^{i,j}}d, j=1,\dots,k\bigg\}\subset C_n.
$$
It is easy to see that 
$$
s_n(x)=s_{n+1}(x)=\cdots\quad\text{whenever $x\notin [F_0*C_1*\cdots *C_{n-1}*C_n^{\max}]_n$.}
\tag 5-1
$$ 
Let $Y_n:=\{(x,t)\in X_0\times\Bbb Z_+\mid 0\le t<s_n(x)\}$.
Then the map 
$$
\phi_n:Y_n\ni ((\boldsymbol f_0,\boldsymbol c_1,\boldsymbol c_2,\dots),t)\mapsto (t\cdot \boldsymbol f_0*\boldsymbol c_1*\cdots *\boldsymbol c_{n-1}, \boldsymbol c_{n},\boldsymbol c_{n+1},\dots)\in X_n
$$
is a homeomorphism of $Y_n$ on $X_n$ for each $n\ge 0$. 
Of course, $Y_0\subset Y_1\subset\cdots$.
We now let 
$$
s:=\sup_{n\ge 0}s_n\quad \text{and}\quad 
Y:=\{(x,t)\in X_0\times\Bbb Z_+\mid 0\le t<s(x)\}.
$$
We call $s$ {\it the spacer map}.
It takes values in $\Bbb Z_+\cup\{+\infty\}$
and it is  lower semicontinuous.
Denote by $\Cal D$ the subset of all $x=(f_0,c_1,c_2,\dots)\in X_0$ such that
 $c_n\in C_n^{\max}$
eventually.
Then $\Cal D$ is an $F_\sigma$-subset of $X_0$ and $\mu(\Cal D)=0$.
In fact, it is easy to verify that $\Cal D$ is countable\footnote{This follows from the fact that there are only finitely many points $(f_0,c_1,c_2,\dots)\in X_0$ with $c_n\in C_n^{\max}$ for each $n>0$.}.
It follows from \thetag{5-1} that the spacer map is continuous when restricted to $X_0\setminus \Cal D$.
The map 
$$
\phi:Y\supset Y_n\ni y\mapsto\phi_n(y)\in X_n\subset X, \qquad n\ge 0,
$$
is a homeomorphism of $Y$ onto $X$.
We define an equivalence relation $\Cal Y$ on $Y$ by setting $(y,i)\sim_{\Cal Y} (y',i')$ if $(y,y')\in\Cal R\cap (X_0\times X_0)$, where $\Cal R$ stands for the $(C,R)$-equivalence relation on $X$.
Then $\phi\times\phi$ maps $\Cal Y$ onto $\Cal R$.
We now define a measure $\nu$ on $Y$ by setting
$d\nu(x,i):=d\mu(x)$, $(x,i)\in Y$.
It is easy to verify that $\mu\circ\phi=\nu$.
Hence  $\mu$ is finite if and only if $\int_{X_0}s\,d\mu<\infty$.
By \thetag{5-1}, if $s(x)=+\infty$  then $x\in\Cal D$.
Hence the spacer map is finite $\mu$-almost everywhere on $X_0$.
We also note that
$$
\phi^{-1}T_1\phi (x,t)=
\cases
(x,t+1)& \text{if $t+1<s(x)$ and}\\
((T_1)_{X_0}x, t)&\text{otherwise}.
\endcases
$$
Thus we showed the following proposition.

\proclaim{Proposition 5.1}
 $T_1$ is a tower transformation built under the spacer map over the base $(T_1)_{X_0}$ which is an ergodic transformation of  exact finite rank. 
\endproclaim

It remains to note that we  consider the transformations of exact rank  as  ``higher rank'' analogues of the transformations with pure point spectrum.

\definition{Definition 5.2} We say that $T$ is of {\it balanced} finite rank it is isomorphic to a $(C,F)$-action such that  $(\alpha)$--$(\gamma)$  are satisfied and  $\Lambda_i>0$ for each $i\in\{1,\dots,k\}$.
\enddefinition

In other words, $T$ is of balanced finite rank  if there is a refining approximating sequence of $T$-castles such that the induced (finite measure preserving) transformation $T_1\restriction(\bigsqcup_{j=1}^kB_0^j)$ is of balanced exact finite rank (see Definition~3.1).

Now we prove a higher rank analogue of Lemma~2.2.

\proclaim{Proposition 5.3}
If $T$ is of balanced finite rank then for each pair of cylinders $[A]_l,[B]_l$ in $X$, we have
$$
\lim_{n\to\infty}\frac{\sum_{m=0}^{n-1}\mu([A]_l\cap T_m[B]_l)}{\sum_{m=0}^{n-1}\mu([F_0]_0\cap T_m[F_0]_0)}=
\mu([A]_l)\mu([B]_l).
$$
\endproclaim

\demo{Proof}
We  first note that for each  pair $i,j\in\{1,\dots,k\}$,
$$
\align
\sum_{m=0}^{n-1}\mu([F_0^i]_l\cap T_m[F_0^j]_l)&=
{\sum_{a,b=1}^k\sum_{m=0}^{n-1}\mu([(0,i)*C^{i,a}_{l+1}]_{l+1}\cap T_m[(0,j)*C^{j,b}_{l+1}]_{l+1})}\\
&=
\sum_{a,b=1}^k\sum_{\boldsymbol c\in C_{l+1}^{i,a}}\sum_{\boldsymbol d\in C_{l+1}^{j,b}}\sum_{m=0}^{n-1}\mu([(0,a)]_{l+1}\cap T_{ d+m- c}[(0,b)]_{l+1}).
\endalign
$$
It is easy to see that for all $a,b\in\{1,\dots,k\}$,  and $s\in\Bbb Z$,
$$
\sum_{m=0}^{n-1}\mu([(0,a)]_{l+1}\cap T_{m+s}[(0,b)]_{l+1})\sim_{n\to\infty}
\sum_{m=0}^{n-1}\mu([(0,a)]_{l+1}\cap T_{m}[(0,b)]_{l+1}).
$$
Since $\# C_{l+1}^{i,a}=r_{l+1}^{i,a}$ and $\{(0,a)\}=F_0^a$ for all $i,a\in\{1,\dots,k\}$, we obtain the following  equivalence
$$
\sum_{m=0}^{n-1}\mu([F_0^i]_l\cap T_m[F_0^j]_l)
\sim_{n\to\infty}
\sum_{a,b=1}^k r_{l+1}^{i,a} r_{l+1}^{j,b}\sum_{m=0}^{n-1}\mu([F_0^a]_{l+1}\cap T_m[F_0^b]_{l+1}).\tag 5-2
$$
Let $C:=F_0*C_1*\cdots *C_{l}$.
Then we have
$$
\align
\frac{\sum_{m=0}^{n-1}\mu([A]_l\cap T_m[B]_l)}{\sum_{m=0}^{n-1}\mu([F_0]_0\cap T_m[F_0]_0)} 
&=
\frac{\sum_{i,j=1}^k\sum_{\boldsymbol a\in A^i,\boldsymbol b\in B^j}\sum_{m=0}^{n-1}\mu([\boldsymbol a]_l\cap T_m[\boldsymbol b]_l)}{\sum_{i,j=1}^k\sum_{\boldsymbol c\in C^i,\boldsymbol d\in C^j}\sum_{m=0}^{n-1}\mu([\boldsymbol c]_l\cap T_m[\boldsymbol d]_l)} \\
&=
\frac{\sum_{i,j=1}^k\sum_{\boldsymbol a\in A^i,\boldsymbol b\in B^j}\sum_{m=0}^{n-1}\mu([F_0^i]_l\cap T_{b+m- a}[F_0^j]_l)}
{\sum_{i,j=1}^k\sum_{\boldsymbol c\in C^i,\boldsymbol d\in C^j}\sum_{m=0}^{n-1}\mu([F_0^i]_l\cap T_{ d+m-c}[F_0^j]_l)} .
\endalign
$$
It follows from this and \thetag{5-2} that there is a limit
$$
\lim_{n\to\infty}\frac{\sum_{m=0}^{n-1}\mu([A]_l\cap T_m[B]_l)}{\sum_{m=0}^{n-1}\mu([F_0]_0\cap T_m[F_0]_0)}=
\frac{\sum_{i,j=1}^k\# A^i\#B^j  \sum_{p,q=1}^kr^{i,p}_{l+1}r^{j,q}_{l+1} }
{\sum_{i,j=1}^k\# C^i\#C^j  \sum_{p,q=1}^k   r^{i,p}_{l+1}r^{j,q}_{l+1}  }.\tag5-3
$$
Since $T$ is of balanced finite rank, 
there is $D>0$ such that
$$
\max_{l\ge 1}\max_{1\le i\le k}\frac{\sum\lambda_l^i}{\lambda_l^j}<D.
$$
It now follows from  this inequality, $(\gamma)$  and \thetag{5-3} that for each $\epsilon>0$, there is $L$ such that if $l>L$ then
$$
\align
\lim_{n\to\infty} \frac{\sum_{m=0}^{n-1}\mu([A]_l\cap T_m[B]_l)}{\sum_{m=0}^{n-1}\mu([F_0]_0\cap T_m[F_0]_0)}&=
\frac{\sum_{i,j=1}^k\# A^i\#B^j  \sum_{p,q=1}^k(\Lambda_{i}\pm \epsilon)(\Lambda_{j}\pm\epsilon)}
{\sum_{i,j=1}^k\# C^i\#C^j   \sum_{p,q=1}^k(\Lambda_{i}\pm \epsilon)(\Lambda_{j}\pm\epsilon) }\\
&=
\frac{\sum_{i,j=1}^k\# A^i\#B^j\lambda_{l}^i\lambda^j_{l}(1\pm 2\epsilon D)^2}{\sum_{i,j=1}^k\# C^i\#C^j\lambda_{l}^i\lambda^j_{l}(1\pm 2\epsilon D)^2}\\
&=\frac{\sum_{i,j=1}^k\mu([A^i]_l)\mu([B^j]_l)(1\pm 2\epsilon D)^2}
{\sum_{i,j=1}^k\mu([C^i]_l)\mu([C^j]_l)(1\pm 2\epsilon D)^2}\\
&=\mu([A]_l)\mu([B]_l)(1\pm 10\epsilon D)^2.
\endalign
$$
Since every cylinder $J$ in $X$ can be presented as $J=[J_l]_l$ for each sufficiently large $l$ and some finite subset $J_l\subset \Bbb Z$, we are done.
\qed
\enddemo

We now let $h_l:=\min_{1\le i\le k}\# F_l^i$ for  each $l\ge 1$.
The following assertion is a higher rank  analogue of Lemma~2.3.

\proclaim{Lemma 5.4}
Under the condition of Proposition~5.3,
for arbitrary subsets $A$ and $B$ in $X_0$, 
$$
\frac{\sum_{m=0}^{ h_l-1}\mu(A\cap T_mB)}{\sum_{m=0}^{ h_l-1}\mu(X_0\cap T_mX_0)}\le
\frac{4k\min(\mu(A),\mu(B))}{\min_{1\le j\le k}\delta_j}.\tag5-4
$$
\endproclaim
\demo{Proof}
Let $C:=F_0*C_1*\cdots *C_l$.
For each subset $J\subset X_0=[F_0]_0$ and  element $\boldsymbol  c\in C$, we set $[\boldsymbol c]_{l,J}:=[\boldsymbol c]_l\cap J$.
Then $J=\bigsqcup_{\boldsymbol c\in C}[\boldsymbol c]_{l,J}$.
We now have 
$$
\align
\sum_{m=0}^{h_l-1}\mu(A\cap T_m B) &\le\sum_{m=0}^{ h_l-1}\sum_{\boldsymbol c,\boldsymbol d\in C}\mu([\boldsymbol c]_{l,A}\cap T_m[\boldsymbol d]_{l})\\
&=
\sum_{\boldsymbol c\in C}\sum_{j=1}^k\sum_{\boldsymbol d\in C^j}\sum_{m=0}^{ h_l-1}\mu([\boldsymbol c]_{l,A}\cap T_{d} [(m,j)]_{l})\\
&\le
\sum_{\boldsymbol c\in C}\sum_{j=1}^k\sum_{\boldsymbol d\in C^j}\mu([\boldsymbol c]_{l,A}\cap T_{ d} [F_l^j]_{l})\\
&\le
k\max_{1\le j\le k}\# C^j\sum_{\boldsymbol c\in C}\mu([\boldsymbol c]_{l,A})\\
&=k\mu(A) \max_{1\le j\le k}\# C^j.
\endalign
$$
In a similar way,
$$
\align
\sum_{m=0}^{h_l-1}\mu(A\cap T_m B) &\le
\sum_{\boldsymbol d\in C}\sum_{j=1}^k
\sum_{\boldsymbol c\in C^j}\sum_{m=0}^{h_l-1}\mu(T_{ c+1-h_l}[(h_l-1-m,j)]_{l}\cap [\boldsymbol d]_{l,B})\\
&\le
\sum_{\boldsymbol d\in C}\sum_{j=1}^k
\sum_{\boldsymbol c\in C^j}\mu(T_{c+1-h_l}[F^j_l]_{l}\cap [\boldsymbol d]_{l,B})\\
&\le
k\max_{1\le j\le k}\# C^j\sum_{\boldsymbol d\in C}\mu([\boldsymbol d]_{l,B})\\
&\le
k\mu(B) \max_{1\le j\le k}\# C^j.
\endalign
$$
Hence
$$
\sum_{m=0}^{h_l-1}\mu(A\cap T_mB)\le k\min(\mu(A),\mu(B))\max_{1\le j\le k}\# C^j .\tag5-5
$$
Choose $j_l$ such that $h_l=\# F_l^{j_l}$, $l\ge 1$.
Then we have
$$
\aligned
\sum_{m=0}^{h_l-1}\mu([F_0]_0\cap T_m [F_0]_0) &\ge
\sum_{i=1}^k\sum_{m=0}^{ h_l-1}\mu([C^i]_{l}\cap T_m[C^i]_{l})\\
&\ge
\sum_{i=1}^k\sum_{c,d\in C^i, \boldkey c\ge\boldkey d}\sum_{m=0}^{ h_l-1}\mu([c]_{l}\cap T_m[d]_{l})\\
&\ge
\sum_{i=1}^k\sum_{c,d\in C^i, \boldkey c\ge\boldkey d}\mu([c]_{l}\cap T_{\boldkey d} [F_l^{j_l}]_{l})
\\
&
\ge
\sum_{i=1}^k\sum_{c,d\in C^i, \boldkey c\ge\boldkey d}\mu([c]_{l})
\\
&=\sum_{i=1}^k\frac{\# C^i+1}2\#C^i\lambda_l^i\\
&\ge\sum_{i=1}^k\frac{\# C^i\mu([C^i]_l)}2.
\endaligned
\tag5-6
$$
We note that $[C^i]_l=X_0\cap[F_l^i]_l$.
Hence $(\beta)$ yields that $\mu([C^i]_l)\ge\delta_i/2$ whenever $l$ is large enough. 
This inequality, \thetag{5-5} and \thetag{5-6} imply \thetag{5-4}.  \qed
\enddemo

We now state the main result of this section.
It is a generalization of Theorem~2.4 to the transformations of balanced finite rank.

\proclaim{Theorem 5.5} Let $T$ be an ergodic i.m.p. $\Bbb Z$-action of balanced finite rank.
Let $h_l:=\min_{1\le i\le k}\# F_l^i$ for  each $l\ge 1$.
Then
for arbitrary Borel subsets $A$ and $B$ in $X_0$, 
$$
\lim_{\l\to\infty}\frac{\sum_{m=0}^{h_l-1}\mu(A\cap T_mB)}{\sum_{m=0}^{h_l-1}\mu(X_0\cap T_mX_0)}=
\mu(A)\mu(B),
$$
i.e. $T$ is subsequence weakly rationally ergodic.
Hence
the ergodic  transformations of balanced finite rank are non-squashable. 
\endproclaim

\demo{Proof}
The theorem is proved in the very same way as Theorem~2.4 but one need to apply Proposition~5.3 and Lemma~5.4 instead of Lemma~2.2 and Lemma 2.3 respectively.
\qed
\enddemo

\head 6. Partial rigidity  for transformations of finite rank
\endhead

In this section we consider only finite measure preserving transformations  ($\Bbb Z$-actions) of finite rank.

We recall that a  measure preserving transformation $S$ of a probability space $(Y,\goth C,\nu)$ is called {\it partially rigid} if there is $\delta>0$ and an increasing  sequence of integers $(m_n)_{n=1}^\infty$ such that
$\mu(A\cap S^{m_n}A)\ge\delta\mu(A)$ for each subset $A\in\goth C$.

The following theorem refines an unpublished result of Rosenthal \cite{Ro}   that the transformations of exact finite rank are not mixing (cf. \cite{Be--So}).

\proclaim{Theorem 6.1}
Let $T=(T_m)_{m\in\Bbb Z}$ be  an ergodic $\Bbb Z$-action of  an exact finite rank.
Then $T$ is partially rigid.
\endproclaim

\demo{Proof\footnote{ Since our original proof  of Theorem 6.1 is rather long we replace it   here with a more elegant  (and slightly modified) Ryzhikov's proof reconstructed from \cite{Ry}.}} 
Let $T$ be of exact rank at most $k$.
Let $(X,\mu)$ stand for  the space of $T$.
We will use the notation from Definition~3.1.
Denote by $h_{n,j}$ the hight of the $j$-th $T$-tower of the $n$-th $T$-castle. 
Thus we have $F_n^j=\{0,1,\dots,h_{n,j}-1\}$.
Fix $n>0$.
Select $j_0\in\{1,\dots,k\}$ such that $\mu(B_n^{j_0})=\max_{1\le j\le k}\mu(B_n^j)$.
Since $T_{h_{n,j_0}}B_n^{j_0}\subset\bigsqcup_{j=1}^kB_n^j$, there is $j_1\in\{1,\dots,k\}$ such that
$$
\mu(B_n^{j_1}\cap T_{h_{n,j_0}}B_n^{j_0})=\max_{1\le j\le k}\mu(B_n^{j}\cap T_{h_{n,j_0}}B_n^{j_0})\ge\frac{\mu(B_n^{j_0})}{k}.
$$
Since $T_{h_{n,j_1}}B_n^{j_1}\subset\bigsqcup_{j=1}^kB_n^j$,
there is $j_2\in\{1,\dots,k\}$ such that
$$
\mu(B_n^{j_2}\cap T_{h_{n,j_1}}(B_n^{j_1}
\cap T_{ h_{n,j_0}}B_n^{j_0}))=\max_{1\le j\le k}\mu(B_n^{j}\cap T_{h_{n,j_1}}(B_n^{j_1}
\cap T_{ h_{n,j_0}}B_n^{j_0}))\ge\frac{\mu(B_n^{j_0})}{k^2}.
$$
Continuing this process $k$ times, we define integers $j_0,\dots,j_k\in\{1,\dots,k\}$.
Hence there are  integers $a,b$ such that $0\le a<b\le k$ and $j_a=j_b$.
Relabeling the towers of the $n$-th $T$-castle,
we may assume without loss of generality that $j_a=1$.
We now let $m_n:=h_{n,j_a}+\cdots+h_{n,j_b-1}$.
Then
$$
\mu(T_{m_n}B_{n}^{1}\cap B_n^{1})\ge\frac{\mu(B_n^{j_0})}{k^k}\ge\frac{\mu(B_n^{1})}{k^k}.
\tag6-1
$$
Let $A\subset X$ be a union of levels of the $n$-th $T$-castle.
Then there is a subset $J\subset F_n^1$ such that $A\cap W_n^1=\bigsqcup_{j\in J}T_jB_n^1$.
Then
$$
\aligned
\mu(T_{m_n}A\cap A)&\ge\mu(T_{m_n}(A\cap W_n^1)\cap (A\cap W_n^1))\\
&\ge\sum_{j\in J}\mu(T_{m_n}(A\cap T_jB_n^1)\cap(A\cap T_jB_n^1) )\\
&=\sum_{j\in J}\mu(T_{m_n}B_n^1\cap B_n^1 ).
\endaligned
\tag6-2
$$
Since $T$ is of exact rank at most $k$, there is $\delta>0$ such that
$\liminf_{n\to\infty}\mu(W_n^1)=\delta$.
It follows from Lemma~3.6 that $\liminf_{n\to\infty}\mu(A\cap W_n^1)=\delta\mu(A)$.
This and \thetag{6-1} with \thetag{6-2} yield that $\liminf_{n\to\infty}\mu(T_{m_n}A\cap A)\ge\frac{\mu(A)}{k^k}$.
It follows from the standard lemma below that $T$ is partially rigid.
\qed

\comment

Suppose first that the corresponding sequence of $T$-castles (see Definition~3.1) refines.
Then we can assume without loss of generality that  $T$ is constructed via a standard cutting-and-stacking  inductive construction in the following way.
For each $n\ge 0$ and  each $i\in\{1,\dots,k\}$, we cut  the $i$-th $T$-tower of the $n$-th castle into finitely many $T$-subtowers.
We call these $T$-subtowers the {\it copies of the $i$-th $T$-tower}.
Then we build  $k$ new $T$-towers by stacking  these $T$-subtowers in such a way that
$\bigsqcup_{i=1}^kW_n^i=\bigsqcup_{i=1}^kW_{n+1}^i$.
We recall that $W_n^i$ denotes the union of all levels in the $i$-th $T$-tower of the $n$-th $T$-castle.
Of course, this is only possible if for each $i\in\{1,\dots,k\}$, the measures of the levels from the $T$-subtowers included into the $i$-th tower of the $(n+1)$-th  $T$-castle are same.
In view of Proposition~3.4 and Corollary~3.7 we may assume that for   all $l,m\in\{1,\dots,k\}$
and each subset $A\subset X$,
there exist limits 
$$
\delta_l:=\lim_{n\to\infty}\mu(W_n^l)>0\quad\text{and}\quad \lim_{n\to\infty}\mu(A\cap W_n^l\cap W_{n+1}^m)=\delta_l\delta_m\mu(A).\tag6-1
$$
We  let $\delta:=\min_{1\le l\le k}\delta_l$.
Fix $j\in\{1,\dots,k\}$.
The number of copies of the $l$-th $T$-tower of the $n$-th $T$-castle inside the $j$-th $T$-tower of the $(n+1)$-th $T$-castle is denoted by $N_n^l$.
For each $n$, we choose $i_n$ so that
$N_n^{i_n}=\max_{1\le l\le k}N_n^l$.
We say that a copy of the $i_n$-th $T$-tower of the $n$-th $T$-castle in the $j$-th $T$-tower of the $(n+1)$-th $T$-castle is {\it good} if there are no more than $k$ copies of $T$-towers of the $n$-th $T$-castle between it and the subsequent (in the natural order on  the $j$-th $T$-tower) copy of $i_n$-th $T$-tower.
Otherwise we say that the copy is {\it bad}.
Consider separately two cases.

(A) Suppose that there is an increasing  sequence $(n_m)_{m\to\infty}$ of positive reals such that the number of good copies of the $i_{n_m}$-th $T$-tower of the $n_m$-th $T$-castle in  the $j$-th $T$-tower of the $(n_m+1)$-th $T$-castle is more than $0.3N_{n_m}^{i_{n_m}}$.
Let $h_l^s$ denote the height (i.e. the number of levels) of the $s$-th $T$-tower in the $n$-th $T$-castle. 
We set
$$
D_m:=\{ h_{n_m}^{j_1},h_{n_m}^{j_1}+h_{n_m}^{j_2},\dots, h_{n_m}^{j_1}+\cdots +h_{n_m}^{j_k}\mid 1\le j_1,\dots,j_k\le k \}.
$$
Since $T$ is of rank at most $k$ without spacers,
for each good copy of the $i_{n_m}$-th $T$-tower, there is $d_m\in D_m$ such that $T_{d_m}$ moves this good copy to the subsequent  copy  of the $i_{n_m}$-th $T$-tower in the $j$-th $T$-tower of the $(n_m+1)$-th $T$-castle.
Let $K:=\sup_{m>0}\# D_m$.
Then $K\le 1+k+\cdots+k^k$.
Therefore there is $g_m\in D_m$ such that $T_{g_m}$ moves at least $\frac {0.3}{K}N_{n_m}^{i_{n_m}}$ of good copies of the $i_{n_m}$-th $T$-tower into copies of the $i_{n_m}$-th $T$-tower. 
We denote by $V_{m}$ the union of these good copies.
Then $V_m\subset W_{n_m}^{i_{n_m}}\cap W_{n_m+1}^j$ and 
$\mu(V_m)\ge \frac{0.3}K\mu(W_{n_m}^{i_{n_m}}\cap W_{n_m+1}^j)$.
Moreover, for each subset $L\subset X$ which is a union of levels in  the $i_{n_m}$-th $T$-tower of the $n_m$-th $T$-castle,
$$
\mu(L\cap V_m)\ge\frac{0.3}K\mu(L\cap W_{n_m}^{i_{n_m}}\cap W_{n_m+1}^j).
\tag6-2
$$
Let $I\subset X$ be  a  level  in some $T$-castle.
Then $J:=I\cap W_{n_m}^{i_{n_m}}$ is a union of levels of the $i_{n_m}$-th $T$-tower of the
$n_m$-th $T$-castle for each  sufficiently large $m$.
It follows from the definition of $V_m$ that $T_{g_m}(J\cap V_m)\subset J$.
From this fact and \thetag{6-2} we deduce that 
$$
\mu(T_{g_m}I\cap I)\ge \mu(T_{g_m}(J\cap V_m)\cap J)=\mu(J\cap V_m)\ge\frac{0.3\mu(I\cap W_{n_m}^{i_{n_m}}\cap W^j_{n_m+1})}{K}.
$$
Applying  \thetag{6-1}  we obtain that $\mu(T_{g_m}I\cap I)\ge  \frac{0.2\delta^2\mu(I)}K$
eventually in $m$.

(B) Suppose that there is an increasing  sequence $(n_m)_{m\to\infty}$ of positive reals such that the number of bad copies of the $i_{n_m}$-th $T$-tower of the $n_m$-th $T$-castle in  the $j$-th $T$-tower of the $(n_m+1)$-th $T$-castle is more than $0.3N_{n_m}^{i_{n_m}}$.
 Among $k+1$ copies of  $T$-towers of the $n$-th $T$-castle following each bad copy, there are two copies of the very same $T$-tower which is different from the $i_{n_m}$-th.
 Therefore there is $p_m\in\{1,\dots,k\}$ and $g_m\in D_m$ such that 
 \roster
 \item"---"
 there exist $\frac {0.3}{kK}N_{n_m}^{i_{n_m}}$ pairs of  copies of the $p_m$-th $T$-tower from the $n_m$-th $T$-castle inside the $j$-th $T$-tower of the $(n_m+1)$-th $T$-castle and
 \item"---"
 $T_{g_m}$ moves the first copy of the $p_m$-th $T$-tower from each such pair to the second copy of  the $p_m$-th $T$-tower from the pair.
\endroster
Denote by $V_m$ the union of the first copies of the $p_m$-th $T$-tower from all these  pairs.
then $V_m\subset W_{n_m}^{p_m}\cap W_{n_m+1}^j$ and $\mu(V_m)\ge\frac{0.3}{kK}\mu(W_{n_m}^{p_m}\cap W_{n_m+1}^j)$ because $N_{n_m}^{i_{n_m}}\ge N_{n_m}^{p_m}$.
Now repeating the reasoning from (A) almost verbally we obtain that $\mu(T_{g_m}I\cap I)> \frac{0.2\delta^2\mu(I)}{kK}$ for all sufficiently large  $m$ for each level $I$ of the $(n_m-1)$-th castle.

From (A), (B) and Lemma~6.2 below we deduce that $T$ is partially rigid.

Consider now the general case, i.e. when  the sequence of $T$-castles from the definition of $T$ does not necessarily refine.
 Then  we can pass to a subsequence that {\it almost} refines, i.e. for each $\epsilon>0$ there is $N>0$ such that every level of the $n$-th $T$-castle coincides with a union of levels in the $(n+1)$-th $T$-castle  up to $\epsilon$ in measure whenever $n>N$.
It is easy to see that the inequalities from the above argument (the refining case) hold in the almost refining case as well.
This yields that $T$ is partially rigid in the general case.

\endcomment
\enddemo

\proclaim{Lemma 6.2} Let $T$ be a probability preserving $\Bbb Z$-action  of funny rank at most $k$.
Let there exist a sequence $(d_n)_{n=1}^\infty$ of positive integers and $\eta>0$ such that
$\liminf_{n\to\infty}\mu(J\cap T_{d_n}J)\ge\eta\mu(J)$ for each level $J$ of the  $l$-th $T$-castle for each $l\ge 0$.
Then $T$ is partially rigid.
\endproclaim

We now generalize the concept of exact finite rank.

\definition{Definition 6.3}
Let $T$ be an ergodic probability preserving $\Bbb Z$-action.
We say that $T$ is  of {\it quasi-exact}  rank at most $k$ if there is an refining sequence of $T$-castles (as in Definition~3.1), $\delta>0$ and  $R>0$ such that
for each $n>0$ and each $j\in\{1,\dots,k\}$, the number of levels (spacers) between two consecutive copies of $T$-towers from the $n$-th $T$-castle in the $j$-th tower of the $(n+1)$-th $T$-castle
is uniformly bounded by $R$ and $\inf_{n\ge 0}\min_{1\le i\le k}\mu(W_n^i)>\delta$.
\enddefinition

Slightly modifying the proof of Theorem~6.1\footnote{Just 
take into account that we now have the inclusion $T_{h_{n,j}}B_n^j\subset\bigcup_{r=0}^RT_{-r}(\bigsqcup_{j=1}^k B_n^j)$ instead of $T_{h_{n,j}}B_n^j\subset \bigsqcup_{j=1}^k B_n^j$, $j=1,\dots,k$.
} we obtain the following theorem.

\proclaim{Theorem 6.4} Let $T$ be an ergodic  $\Bbb Z$-action of quasi-exact finite rank.
 Then $T$ is partially rigid. 
\endproclaim

We now adapt the  definition of consecutive ordering from the theory of Brattelli-Vershik systems  (see \cite{Du}) to the context of measurable systems of finite rank.

\definition{Definition 6.5} Let $T$ be a $\Bbb Z$-action of  rank at most $k$  without spacers and the corresponding sequence of  $T$-castles (see Definition~3.1) refines.
We say that $T$ {\it satisfies the CO-condition} if
given arbitrary $n\ge 0$ and $i,j,l\in\{1,\dots,k\}$, if a copy of the $l$-th $T$-tower from the $n$-th  $T$-castle is between 2 copies of the $i$-th $T$-tower inside the $j$-th $T$-tower of the $(n+1)$-th $T$-castle then $l=i$.
\enddefinition

We note  that the Bratteli-Vershik maps  corresponding to the minimal IETs satisfy the CO-condition \cite{Gj--Jo}.

It was shown in \cite{Be--So} that if $T$ satisfies  the CO-condition and  an additional ``non-degeneracy'' condition (as in Theorem~6.6 below) then $T$ is not mixing.
We prove a stronger result.

\proclaim{Theorem 6.6} Let $T$ be an ergodic $\Bbb Z$-action of rank at most $k$ without spacers.
Let  $T$ satisfy the CO-condition.
Suppose also that 
if  for some $n\ge 0$ and $i,j\in\{1,\dots,k\}$, if a copy of the  $i$-th $T$-tower from the $n$-th $T$-castle is contained  inside the $j$-th $T$-tower from the $(n+1)$-th $T$-castle then at least one more copy of the $i$-th $T$-tower is contained inside the $j$-th $T$-tower.
 Then $T$ is partially rigid.
\endproclaim

\demo{Proof}
Let $(X,\mu)$ stand for the space of $T$.
Since $\sum_{i=1}^k\mu(W_n^i)=1$ for each $n\ge 0$, there is a subsequence $(n_m)_{m=1}^\infty$ of positive integers and an integer $p\in\{1,\dots,k\}$ such that
$\lim_{m\to\infty}\mu(W_{n_m}^p)=\delta>1/k$.
For each $l\ge 0$ and each level $I$ of the $l$-th $T$-castle, the intersection of $I$ with the $p$-th tower of the $n_m$-th $T$-castle is the union of some levels of this tower for all sufficiently large $m$.
Let $h_{m}$ stand for the height of the $p$-th tower of the $n_m$-th $T$-castle.
Passing to a further subsequence we may assume withal loss of generality that
there is $q\in\{1,\dots,k\}$ such that $\mu(W_{n_m}^p\cap W_{n_m+1}^q)\ge \mu(W_{n_m}^p)/k$ for all $m$.
Denote by $V_m$ the union of all but  the top  one copies of the $p$-th tower of the $n_m$-th $T$-castle in the $q$-th tower of the $(n_m+1)$-th $T$-castle.
It follows from the condition of the theorem that $T^{h_m}$ moves each copy  of the $p$-th $T$-tower in $V_m$
onto the adjacent  (from above) copy of the same $T$-tower inside the $q$-th $T$-tower from $(n_m+1)$-th $T$-castle.
Therefore $ T_{h_m}(I\cap W_{n_m}^p\cap V_m)\subset I\cap W_{n_m}^p$.
Since there are no less then 2 such copies inside the $q$-th $T$-tower,  we obtain that
$$
\align
\mu( T_{h_m}I\cap I) &\ge\mu( T_{h_m}(I\cap W_{n_m}^p\cap V_m) \cap I)\\
&={\mu(I\cap W_{n_m}^p\cap V_m)}\\
&\ge\frac 12 \mu(I\cap W_{n_m}^p\cap W_{n_m+1}^q)\\
&\ge\frac 1{2k} \mu(I\cap W_{n_m}^p).
\endalign
$$
It now follows from Lemma~3.6 that
$\mu( T_{h_m}I\cap I)\ge\frac\delta{3k}\mu(I)$ eventually in $m$.
It remains to apply Lemma~6.2.
\qed
\enddemo

\remark{Remark \rom{6.7}} It is possible to generalize Theorem~6.6 (with a slight only modification of the proof) in the following way: drop the assumption that $T$ is constructed without spacers and replace
 the CO-condition in the statement of the theorem with the following two conditions.
\roster
\item"CO$_1$"  There is $R>0$ such that
the number of levels (spacers) between two neighboring copies of $T$-towers from the $n$-th $T$-castle in every tower of the $(n+1)$-th $T$-castle
is uniformly bounded by $R$.
\item"CO$_2$" There is $L>0$ such that
the number of copies of $T$-towers from the $n$-th $T$ castle between every  two
consecutive copies of the $i$-th $T$-tower 
 inside each $T$-tower of the $(n+1)$-th $T$-castle, $i=1,\dots,k, $ is uniformly bounded by $L$.
\endroster
We leave details of the proof to the reader.
\endremark

 A.~Katok proved  in \cite{Ka} that the ergodic IETs  are not mixing.
We now show how to deduce  from that proof (or, rather a slight modification of that proof from \cite{KSF}) that they are partially rigid.\footnote{We can not apply Theorem~6.1 to the ergodic IETs because we do not know whether they are of exact finite rank.}

\proclaim{Proposition 6.8}
Let $T$ be an ergodic IET.
Then $T$ is partially rigid.
\endproclaim
\demo{Proof\footnote{See also \cite{Ry} for an alternative proof.}}
Indeed, by \cite{KSF, Chapter, \S\,3, Lemma~1}, there are integers $k>0$, $r_j^{(n)}> 0$, $j=1,\dots,k$, and measurable partitions $(A_1^{(n)},\dots, A_k^{(n)})$ of $X$,   
 $n\in\Bbb N$, such that 
 \roster
 \item"$(i)$"
 $\min_{1\le j\le k}r_j^{(n)}\to\infty$ as $n\to\infty$,
\item"$(ii)$" 
$\text{Leb}(B\triangle(\bigcup_{j=1}^kT^{r_j^{(n)}}( A^{(n)}_j\cap B))\to 0$ as $n\to\infty$
for each Borel subset $B\subset [0,1)$ and
\item"$(iii)$"
$\text{Leb}(TA^{(n)}_j\triangle A^{(n)}_j)\to 0$ as $n\to\infty$ for each $j=1,\dots,k$.
\endroster
Passing to a subsequence, we may assume without loss of generality that
there exist integers $j_n\in\{1,\dots,k\}$ and a limit
\roster
\item"$(iv)$" 
$\lim_{n\to\infty}\text{Leb}(A_{j_n}^{(n)})=:\delta\ge 1/k$.
\endroster 
It follows from~$(ii)$ that $\sum_{j=1}^k\int_{A_{j}^{(n)}}(U_T^{r_j^{(n)}}1_B-1_B)^2\,dx\to 0$ as $n\to\infty$.
Therefore
 $$
 \int_{A_{j_n}^{(n)}}(U_T^{r_{j_n}^{(n)}}1_B-1_B)^2\,dx=\int_{A_{j_n}^{(n)}}(U_T^{r_{j_n}^{(n)}}1_B-2U_T^{r_{j_n}^{(n)}}1_B1_B+1_B)\,dx\to 0
 $$
 as $n\to\infty$.
 Applying $(i)$, $(iii)$, $(iv)$ and Lemma~3.6 we obtain that
$$
\int_{A_{j_n}^{(n)}}(2U_T^{r_{j_n}^{(n)}}1_B1_B-U_T^{r_{j_n}^{(n)}}1_B)\,dx\to \delta\mu(B)\tag6-3
$$
as $n\to\infty$.
Since
$$
\int_{[0,1)}2U_T^{r_{j_n}^{(n)}}1_B1_B\,dx\ge\int_{A_{j_n}^{(n)}}2U_T^{r_{j_n}^{(n)}}1_B1_B\,dx\ge \int_{A_{j_n}^{(n)}}(2U_T^{r_{j_n}^{(n)}}1_B1_B-U_T^{r_{j_n}^{(n)}}1_B)\,dx,
$$
we deduce from \thetag{6-3} that  $\liminf_{n\to\infty}\mu(T^{r_{j_n}^{(n)}}B\cap B)\ge\frac \delta{2}\mu(B)\ge\frac 1{2k}\mu(B)$.
\qed
\enddemo

 \head 7. Open problems
\endhead
\roster
\item
Given an infinite countable amenable discrete  group $G$, is there a finite measure preserving free  action of $G$ which is of funny rank one?
We note that this question is related closely to a basic   problem in the theory of amenable groups: whether every amenable group has a F{\o}lner sequence consisting of monotiles \cite{We2}?
\item
Whether each funny rank-one i.m.p. action of an Abe\-lian  countable discrete group  $G$ is weakly rationally ergodic? The question is especially interesting in the case where $G=\Bbb Z$.
\item 
Whether the squashability of ergodic i.m.p. actions of Abelian groups $G$ is a spectral property?
Is there en ergodic squashable i.m.p. action of $G$ whose spectrum is of finite multiplicity? 
\item
Given an amenable group $G$, consider two classes of  actions of $G$: the class of all  possible i.m.p. $(C,F)$-actions of funny rank one and the class of  i.m.p. $(C,F)$-actions of funny rank one such that the corresponding sequences $(F_n)_{n\ge 0}$ are F{\o}lner.
Do these classes coincide? 
\item
Let  $\goth F$ be a factor (i.e. an invariant $\sigma$-finite $\sigma$-subalgebra) of an ergodic i.m.p. transformation $T$.
Let  $T\restriction \goth F$ be
  non-squashable.
  Is $T$  non-squashable too?
  \item 
  Let  $\goth F$ be a factor  of an ergodic i.m.p. transformation $T$
  and let  $T\restriction \goth F$ be (subsequence) weakly rationally ergodic.
  Is $T$ (subsequence) weakly rationally ergodic?
  \item
  Given an ergodic  transformation of exact finite rank at least 2, is it possible to find a  refining sequence of approximating castles such that each castles fills the entire space (i.e. the union of all levels in the castle equals the entire space)? Of course, the answer is affirmative for the transformations of exact rank one.
  \item
  Give  examples of  ergodic transformations which are of  exact finite rank but not of balanced exact finite rank.
  \item
  Are there ergodic i.m.p. transformations of finite rank but not of  balanced finite rank?
  If yes, are they subsequence weakly rationally ergodic?
\item 
 Are there ergodic  IETs  which are not of exact finite rank? 
   \endroster

\enddocument